%% file: _English.tex
\documentclass{amsart}
\input{Commands}

\begin{document}
\title{Basis Manifold}


\pdfbookmark[1]{Basis Manifold}{TitleEnglish}
\begin{abstract}
In this paper we study a symmetry group of vector space.
Basis manifold is a homogeneous space of a symmetry group.
This concept leads us to the definition
of active and passive transformations on basis manifold.
Active transformation can be expressed as a transformation of
vector space.
Passive transformation gives ability to define
concepts of invariance and of geometrical object.
\end{abstract}

\maketitle

\input{Representation.English}
\input{Basis.English}
\input{Biblio.English}
\input{Index.English}
\input{Symbol.English}
\end{document}

%% file: Commands.tex
\def\Defined{}
\usepackage{footmisc}
\usepackage[all]{xy}
\usepackage{color}
\definecolor{UrlColor}{rgb}{.9,0,.3}
\definecolor{SymbColor}{rgb}{.4,0,.9}
\definecolor{IndexColor}{rgb}{1,.3,.6}
\definecolor{eml1}{rgb}{.8,.1,.1}
\definecolor{eml2}{rgb}{.1,.6,.6}

\usepackage{xr-hyper}
\usepackage[unicode]{hyperref}
\hypersetup{pdfdisplaydoctitle=true}
\hypersetup{colorlinks}
\hypersetup{citecolor=UrlColor}
\hypersetup{urlcolor=UrlColor}
\hypersetup{pdffitwindow=true}
\hypersetup{pdfnewwindow=true}
\hypersetup{pdfstartview={FitH}}

\def\hyph{\penalty0\hskip0pt\relax-\penalty0\hskip0pt\relax}
\def\Hyph{-\penalty0\hskip0pt\relax}

\newcommand{\Basis}[1]{\overline{\overline{#1}}{}}
\newcommand{\Vector}[1]{\overline{#1}{}}
\newcommand{\gi}[1]{\boldsymbol{\textcolor{IndexColor}{#1}}}
\makeatletter
\newcommand{\NameDef}[1]{%
	\expandafter\gdef\csname #1\endcsname%
}%
\newcommand{\ShowSymbol}[1]{%
	\@nameuse{ViewSymbol#1}%
}%
\newcommand{\symb}[3]{%
	\@ifundefined{ViewSymbol#3}{%
		\NameDef{ViewSymbol#3}{\textcolor{SymbColor}{#1}}%
		\NameDef{RefSymbol#3}{\pageref{symbol: #3}}%
		\@namedef{LabelSymbol#3}{\label{symbol: #3}}%
	}{%
		\NameDef{RefSymbol#3}{}%
		\@namedef{LabelSymbol#3}{}%
	}%
	\ifcase#2
	\or
		$\@nameuse{ViewSymbol#3}$%
	\or
		\[\@nameuse{ViewSymbol#3}\]%
	\else%
	\fi%
	\@nameuse{LabelSymbol#3}%
}%
\makeatother

\newcommand{\subs}{${}_*$\Hyph}
\newcommand{\sups}{${}^*$\Hyph}

\newcommand{\CRstar}{{}_*{}^*}
\newcommand{\RCstar}{{}^*{}_*}

\newcommand{\RC}{$\RCstar$\Hyph}
\newcommand{\CR}{$\CRstar$\Hyph}
\newcommand{\drc}{$D\RCstar$\Hyph}
\newcommand{\dcr}{$D\CRstar$\hyph}
\newcommand{\rcd}{$\RCstar D$\Hyph}
\newcommand{\crd}{$\CRstar D$\Hyph}

\newcommand\sT{$\star T$\Hyph}%
\newcommand\Ts{$T\star$\Hyph}%

\renewcommand{\uppercasenonmath}[1]{}

\makeatletter
\newcommand\@dotsep{4.5}
\def\@tocline#1#2#3#4#5#6#7
{\relax
		\par \addpenalty\@secpenalty\addvspace{#2}%
		\begingroup \hyphenpenalty\@M
		\@ifempty{#4}{%
			\@tempdima\csname r@tocindent\number#1\endcsname\relax
		}{%
			\@tempdima#4\relax
		}%
		\parindent\z@ \leftskip#3\relax \advance\leftskip\@tempdima\relax
		\rightskip\@pnumwidth plus1em \parfillskip-\@pnumwidth
		#5\leavevmode\hskip-\@tempdima #6\relax
		\leaders\hbox{$\m@th
		\mkern \@dotsep mu\hbox{.}\mkern \@dotsep mu$}\hfill
		\hbox to\@pnumwidth{\@tocpagenum{#7}}\par
		\nobreak
		\endgroup
}
\makeatother 

\ifx\PrintBook\undefined
	\usepackage{fancyhdr}
	\makeatletter
	\renewcommand{\@indextitlestyle}{%
		\twocolumn[\section{\indexname}]%
		\def\IndexSpace{off}%
	}
	\makeatother 
	\thanks{\href{mailto:Aleks\_Kleyn@MailAPS.org}{Aleks\_Kleyn@MailAPS.org}}
\else
	\makeatletter
	\renewcommand{\@indextitlestyle}{%
		\twocolumn[\chapter{\indexname}]%
		\def\IndexSpace{off}%
		\let\@secnumber\@empty
		\chaptermark{\indexname}%
	}
	\makeatother 
	\email{\href{mailto:Aleks\_Kleyn@MailAPS.org}{Aleks\_Kleyn@MailAPS.org}}
\fi

\ifx\SelectlEnglish\undefined
	\ifx\UseRussian\undefined
		\def\SelectlEnglish{}
	\fi
\fi

\ifx\SelectlEnglish\undefined
	\newcommand\CurrentLanguage{Russian.}%
	\author{Александр Клейн}
	\newtheorem{theorem}{Теорема}[section]
	\newtheorem{corollary}[theorem]{Следствие}
	\theoremstyle{definition}
	\newtheorem{definition}[theorem]{Определение}
	\newtheorem{example}[theorem]{Пример}
	\newtheorem{xca}[theorem]{Exercise}
	\theoremstyle{remark}
	\newtheorem{remark}[theorem]{Замечание}
	
	\newcommand\Gbasis{$G$\Hyph базис}
	\newcommand\Gcoords{$G$\Hyph координат}
	\newcommand\Gspace{$G$\Hyph пространств}
	
	\newcommand\xRef[2]%
	{%
		\ifx\PrintBook\undefined%
		\citeBib{#1}-\ref{#1-Russian-#2}%
		\else%
		\ref{#2}%
		\fi%
	}%
	\newcommand\xEqRef[2]%
	{%
		\ifx\PrintBook\undefined%
		\citeBib{#1}-\eqref{#1-Russian-#2}%
		\else%
		\eqref{#2}%
		\fi%
	}%
	\ifx\PrintBook\undefined
		\newcommand{\BibTitle}{%
			\section{Список литературы}%
		}
	\else
		\newcommand{\BibTitle}{%
			\chapter{Список литературы}%
		}
	\fi
\else
	\newcommand\CurrentLanguage{English.}%
	\author{Aleks Kleyn}
	\newtheorem{theorem}{Theorem}[section]
	
	\theoremstyle{definition}
	\newtheorem{definition}[theorem]{Definition}
	\newtheorem{example}[theorem]{Example}
	
	\theoremstyle{remark}
	\newtheorem{remark}[theorem]{Remark}
	\newcommand\Gbasis{$G$\Hyph basis}
	\newcommand\Gcoords{$G$\Hyph coordinates}
	\newcommand\Gspace{$G$\Hyph space}
	
	\newcommand\xRef[2]%
	{%
		\ifx\PrintBook\undefined%
		\citeBib{#1}-\ref{#1-English-#2}%
		\else%
		\ref{#2}%
		\fi%
	}%
	\newcommand\xEqRef[2]%
	{%
		\ifx\PrintBook\undefined%
		\citeBib{#1}-\eqref{#1-English-#2}%
		\else%
		\eqref{#2}%
		\fi%
	}%
	\ifx\PrintBook\undefined
		\newcommand{\BibTitle}{%
			\section{References}%
		}
	\else
		\newcommand{\BibTitle}{%
			\chapter{References}%
		}
	\fi
\fi

\ifx\PrintBook\undefined
\else
	\numberwithin{section}{chapter}
\fi

\numberwithin{equation}{section}
\numberwithin{figure}{section}
\numberwithin{table}{section}
\numberwithin{Item}{section}
\numberwithin{Hfootnote}{section}

\makeatletter
\newcommand\org@maketitle{}
\let\org@maketitle\maketitle
\def\maketitle{%
	\hypersetup{pdftitle={\@title}}%
	\hypersetup{pdfauthor={\authors}}%
	\hypersetup{pdfsubject=\@keywords}%
	\org@maketitle
}
\def\make@stripped@name#1{%
	\begingroup
		\escapechar\m@ne
		\global\let\newname\@empty
		\protected@edef\Hy@tempa{\CurrentLanguage #1}%
		\edef\@tempb{%
			\noexpand\@tfor\noexpand\Hy@tempa:=%
			\expandafter\strip@prefix\meaning\Hy@tempa
		}%
		\@tempb\do{%
			\if\Hy@tempa\else
				\if\Hy@tempa\else
					\xdef\newname{\newname\Hy@tempa}%
				\fi
			\fi
		}%
	\endgroup
}%
\newenvironment{enumBib}{%
	\BibTitle
	\advance\@enumdepth \@ne
	\edef\@enumctr{enum\romannumeral\the\@enumdepth}\list
	{\csname biblabel\@enumctr\endcsname}{\usecounter
	{\@enumctr}\def\makelabel##1{\hss\llap{\upshape##1}}}
}{%
	\endlist
}

\def\Chapters#1{\ChapterList#1,LastChapter,}%
\def\LastChapter{LastChapter}%
\def\ChapterList#1,{\def\temp{#1}%
	\ifx\temp\LastChapter
	\else
		\@ifundefined{#1}{%
		}{%
			\def\Semafor{on}
		}
		\expandafter\ChapterList
	\fi
}%
\newcommand{\BiblioItem}[3]
{
	\def\Semafor{off}
	\Chapters{#1}
	\ifx\Semafor\ValueOn
		\ifx\IndexState\ValueOff
			\begin{enumBib}
			\def\IndexState{on}
		\fi
		\item \label{bibitem: #2}#3%
	\fi
}
\newcommand{\OpenBiblio}
{
	\def\IndexState{off}
}
\newcommand{\CloseBiblio}
{
	\ifx\IndexState\ValueOn
		\end{enumBib}
		\def\IndexState{off}
	\fi
}

\def\StartCite{[}%
\def\citeBib#1{\showCiteBib#1,endCite,}%
\def\endCite{endCite}%
\def\showCiteBib#1,{\def\temp{#1}%
\ifx\temp\endCite
]%
\def\StartCite{[}%
\else
	\StartCite\ref{bibitem: #1}%
	\def\StartCite{, }%
\expandafter\showCiteBib%
\fi}%

\makeatother 
\newcommand{\arp}{\ar @{-->}}
\newcommand{\bundle}[4]
{
	\def\tempa{}%
	\def\tempb{#3}%
	\def\tempc{#1}%
	\ifx\tempa\tempb
		\ifx\tempa\tempc
			#2
		\else
			\xymatrix{#2:#1\arp[r]&#4}
		\fi
	\else
		\ifx\tempa\tempc
			#2[#3]
		\else
			\xymatrix{#2[#3]:#1\arp[r]&#4}
		\fi
	\fi
}
\newcommand{\AddIndex}[2]%
{%
	{\bf #1}%
	\label{index: #2}%
}%
\newcommand{\Index}[3]%
{%
	\def\Semafor{off}%
	\Chapters{#1}%
	\ifx\Semafor\ValueOn%
		\def\tempa{}%
		\def\tempb{#3}%
		\ifx\IndexState\ValueOff%
			\begin{theindex}%
			\def\IndexState{on}%
		\fi%
		\ifx\IndexSpace\ValueOn%
			\indexspace%
			\def\IndexSpace{off}%
		\fi%
		\item #2%
		\ifx\tempa\tempb%
		\else%
			\pageref{index: #3}%
		\fi%
	\fi%
}%
\newcommand{\SubIndex}[3]
{
	\def\Semafor{off}
	\Chapters{#1}
	\ifx\Semafor\ValueOn
		\subitem #2 \pageref{index: #3}
	\fi
}%

\makeatletter
\newcommand{\Symb}[3]
{
	\def\Semafor{off}
	\Chapters{#1}
	\ifx\Semafor\ValueOn
		\ifx\IndexState\ValueOff
			\begin{theindex}
			\def\IndexState{on}
		\fi
		\ifx\IndexSpace\ValueOn
			\indexspace
			\def\IndexSpace{off}
		\fi
		\item $\@nameuse{ViewSymbol#3}$\ \ #2
		\@nameuse{RefSymbol#3}%
	\fi
}
\def\CiteBibNote{\xdef\@thefnmark{}\@footnotetext}
\makeatother

\newcommand{\SetIndexSpace}%
{%
	\def\IndexSpace{on}%
}%
\def\ValueOff{off}
\def\ValueOn{on}

\newcommand{\OpenIndex}
{
	\def\IndexState{off}
}
\newcommand{\CloseIndex}
{
	\ifx\IndexState\ValueOn
		\end{theindex}
		\def\IndexState{off}
	\fi
}

\def\LastMemo{LastMemo}%
\def\MemoList#1//{\def\temp{#1}%
	\ifx\temp\LastMemo
	\else%
		\par\setlength{\parindent}{12pt}\textcolor{blue}{#1}%
		\expandafter\MemoList%
	\fi%
}%

%% file: Representation.English.tex
\def\texRepresentation{}
\ifx\PrintBook\undefined

\section{Introduction}

This paper was written under the great influence of the book \citeBib{Rashevsky}.
The studying of a homogenous space of a group of symmetry of a vector space
leads us to the definition of a basis of this space
and a basis manifold.
We introduce two types of transformation of a basis manifold:
active and passive transformations. The difference between them is
that the active transformation can be expressed as a transformation of
an original space.
As it is shown in \citeBib{Rashevsky} passive transformation
gives ability to define
concepts of invariance and of geometrical object.

Two opposite points of view about a geometrical object
meet in definition \ref{definition: geometrical object}.
On the one hand we determine coordinates of the geometrical object
relative to a given basis and introduce the law of transition of
coordinates during transformation of the basis.
At the same time we study the set of coordinates
of the geometrical object relative to different bases
as a single whole. This gives us an opportunity
to study the geometrical object without using coordinates.
\else
				\chapter{Representation of Group}
				\label{chapter: Representation of Group}
\fi
			\section{Representation of Group}
			\label{section: Representation of Group}

	\begin{definition}
	\label{definition: nonsingular transformation} 
We call the map
\[
t:M\rightarrow M
\]
\AddIndex{nonsingular transformation}{nonsingular transformation},
if there exists inverse map.
	 \qed
	 \end{definition}

	 \begin{definition}
Transformations is
\AddIndex{left-side transformations}{left-side transformation}
if it acts from left
\[
u'=t u
\]
We denote
\symb{l(M)}1{set of left-side nonsingular transformations}
the set of left-side nonsingular transformations of set $M$.
	 \qed
	 \end{definition}

	 \begin{definition}
Transformations is
\AddIndex{right-side transformations}{right-side transformation}
if it acts from right
\[
u'= ut
\]
We denote
\symb{r(M)}1{set of right-side nonsingular transformations}
the set of right-side nonsingular transformations of set $M$.
	 \qed
	 \end{definition}

We denote
\symb{\delta}1{identical transformation}
identical transformation.

		\begin{definition}
		\label{definition: left-side representation of group}
Let $G$ be group.
We call map
	\begin{equation}
f:G\rightarrow l(M)
	\label{eq: left-side representation of group}
	\end{equation}
\AddIndex{left-side representation of group}
{left-side representation of group}
$G$ in set $M$ if map $f$ holds
	\begin{equation}
f(ab)u=f(a)(f(b)u)
	\label{eq: left-side product of transformations}
	\end{equation}
	\begin{equation}
f(e)=\delta
	\label{eq: left-side identical transformation}
	\end{equation}
		 \qed
		 \end{definition}

		\begin{definition}
		\label{definition: right-side representation of group} 
Let $G$ be group.
We call map
	\begin{equation}
f:G\rightarrow r(M)
	\label{eq: right-side representation of group}
	\end{equation}
\AddIndex{right-side representation of group}
{right-side representation of group}
$G$ in set $M$ if map $f$ holds
	\begin{equation}
uf(ab)=(uf(a))f(b)
	\label{eq: right-side product of transformations}
	\end{equation}
	\begin{equation}
f(e)=\delta
	\label{eq: right-side identical transformation}
	\end{equation}
		 \qed
		 \end{definition}

Any statement which holds for left-side representation of group
holds also for right-side representation.
For this reason we use the common term
\AddIndex{representation of group}{representation of group}
and use notation for left-side representation
in case when it does not lead to misunderstanding.

			\begin{theorem}
			\label{theorem: inverse transformation}
For any $g\in G$
	\begin{equation}
	\label{eq: inverse transformation}
f(g^{-1})=f(g)^{-1}
	\end{equation}
		\end{theorem}
		\begin{proof}
Since \eqref{eq: left-side product of transformations}
and \eqref{eq: left-side identical transformation},
we have
\[
u=\delta u=f(gg^{-1})u=f(g)(f(g^{-1})u)
\]
This completes the proof.
		\end{proof}

		 \begin{theorem}
		 \label{theorem: left-side covariant representation}
Let $l(M)$ be a group with respect to multiplication
	\begin{equation}
(t_1t_2)u=t_1(t_2u)
	\label{eq: left-side covariant product of transformations}
	\end{equation}
and $\delta$ be unit of group $l(M)$.
Let map \eqref{eq: left-side representation of group}
be a homomorphism of group
	\begin{equation}
f(ab)=f(a)f(b)
	\label{eq: left-side covariant representation}
	\end{equation}
Then this map is representation of group $G$ which we call
\AddIndex{left-side covariant representation}
{left-side covariant representation}.
		\end{theorem}
		\begin{proof}
Since $f$ is homomorphism of group, we have $f(e)=\delta$.

Since
\eqref{eq: left-side covariant product of transformations}
and \eqref{eq: left-side covariant representation}, we have
\[
f(ab)u=(f(a)f(b))u=f(a)(f(b)u)
\]

According definition \ref{definition: left-side representation of group}
$f$ is representation.
		\end{proof}

		 \begin{theorem}
		 \label{theorem: right-side covariant representation}
Let $r(M)$ be a group with respect to multiplication
	\begin{equation}
u(t_1t_2)=(ut_1)t_2
	\label{eq: right-side covariant product of transformations}
	\end{equation}
and $\delta$ be unit of group $r(M)$.
Let map \eqref{eq: right-side representation of group}
be a homomorphism of group
	\begin{equation}
f(ab)=f(a)f(b)
	\label{eq: right-side covariant representation}
	\end{equation}
Then this map is representation of group $G$ which we call
\AddIndex{right-side covariant representation}
{right-side covariant representation}.
		\end{theorem}
		\begin{proof}
Since $f$ is homomorphism of group, we have $f(e)=\delta$.

Since
\eqref{eq: right-side covariant product of transformations}
and \eqref{eq: right-side covariant representation}, we have
\[
uf(ab)=u(f(a)f(b))=(uf(a))f(b)
\]

According definition \ref{definition: right-side representation of group}
$f$ is representation.
		\end{proof}

	 \begin{theorem}
	 \label{theorem: left-side contravariant representation}
Let $l(M)$ be a group with respect to multiplication
	\begin{equation}
(t_2t_1)u=t_1(t_2u)
	\label{eq: left-side contravariant product of transformations}
	\end{equation}
and $\delta$ be unit of group $l(M)$.
Let map \eqref{eq: left-side representation of group}
be an antihomomorphism of group
	\begin{equation}
f(ba)=f(a)f(b)
	\label{eq: left-side contravariant representation}
	\end{equation}
Then this map is representation of group $G$ which we call
\AddIndex{left-side contravariant representation}
{left-side contravariant representation}.
		\end{theorem}
		\begin{proof}
Since $f$ is antihomomorphism of group, we have $f(e)=\delta$.

Since
\eqref{eq: left-side contravariant product of transformations}
and \eqref{eq: left-side contravariant representation}, we have
\[
f(ab)u=(f(b)f(a))u=f(a)(f(b)u)
\]

According definition \ref{definition: left-side representation of group}
$f$ is representation.
		\end{proof}

	 \begin{theorem}
	 \label{theorem: right-side contravariant representation}
Let $r(M)$ be a group with respect to multiplication
	\begin{equation}
u(t_2t_1)=(ut_1)t_2
	\label{eq: right-side contravariant product of transformations}
	\end{equation}
and $\delta$ be unit of group $r(M)$.
Let map \eqref{eq: right-side representation of group}
be an antihomomorphism of group
	\begin{equation}
f(ba)=f(a)f(b)
	\label{eq: right-side contravariant representation}
	\end{equation}
Then this map is representation of group $G$ which we call
\AddIndex{right-side contravariant representation}
{right-side contravariant representation}.
		\end{theorem}
		\begin{proof}
Since $f$ is antihomomorphism of group, we have $f(e)=\delta$.

Since
\eqref{eq: right-side contravariant product of transformations}
and \eqref{eq: right-side contravariant representation}, we have
\[
uf(ab)=u(f(b)f(a))=(uf(a))f(b)
\]

According definition \ref{definition: right-side representation of group}
$f$ is representation.
		\end{proof}

	\begin{example}
The group operation determines two different representations on the group:
the \AddIndex{left shift}{left shift, group} which we introduce by the equation
\symb{L(a)b}0{left shift}
	\begin{equation}
b'=\ShowSymbol{left shift}=ab
	\label{eq: left shift}
	\end{equation}
and the \AddIndex{right shift}{right shift, group} which we introduce by the equation
\symb{R(a)b}0{right shift}
	\begin{equation}
b'=\ShowSymbol{right shift}=ba
	\label{eq: right shift}
	\end{equation}
	 \qed
	 \end{example}

		\begin{theorem}
Let representation
	\[
u'=f(a)u
	\]
be covariant representation. Then representation
	\[
u'=h(a)u=f(a^{-1})u
	\]
is contravariant representation.
		\end{theorem}
		\begin{proof}
Statement follows from chain of equations
\[h(ab)=f((ab)^{-1})=f(b^{-1}a^{-1})=f(b^{-1})f(a^{-1})=h(b)h(a)\]
		\end{proof}

	 \begin{definition}
Let $f$ be representation of the group $G$ in set $M$.
For any $v\in M$ we define
\AddIndex{orbit of representation of the group $G$}
{orbit of representation of group} as set
\symb{\mathcal O(v,g\in G,f(g)v)=\{w=f(g)v:g\in G\}}2
{orbit of representation of group}
	 \qed
	 \end{definition}

Since $f(e)=\delta$ we have
$v\in\mathcal{O}(v,g\in G,f(g)v)$.

		\begin{theorem}
		\label{theorem: proper definition of orbit}
Suppose
	\begin{equation}
	\label{eq: orbit, proposition}
v\in\mathcal{O}(u,g\in G,f(g)u)
	\end{equation}
Then
\[
\mathcal{O}(u,g\in G,f(g)u)=\mathcal{O}(v,g\in G,f(g)v)
\]
		\end{theorem}
		\begin{proof}
From \eqref{eq: orbit, proposition} it follows
that there exists $a\in G$ such that
	\begin{equation}
v=f(a)u
	\label{eq: orbit, 1}
	\end{equation}
Suppose $w\in\mathcal{O}(v,g\in G,f(g)v)$. Then
there exists $b\in G$ such that
	\begin{equation}
w=f(b)v
	\label{eq: orbit, 2}
	\end{equation}
If we substitute \eqref{eq: orbit, 1}
into \eqref{eq: orbit, 2} we get
	\begin{equation}
w=f(b)(f(a)u)
	\label{eq: orbit, 3}
	\end{equation}
Since \eqref{eq: left-side product of transformations},
we see that
from \eqref{eq: orbit, 3} it follows
that $w\in\mathcal{O}(u,g\in G,f(g)u)$.
Thus
\[
\mathcal{O}(v,g\in G,f(g)v)\subseteq\mathcal{O}(u,g\in G,f(g)u)
\]

Since \eqref{eq: inverse transformation},
we see that
from \eqref{eq: orbit, 1} it follows
that
	\begin{equation}
u=f(a)^{-1}v=f(a^{-1})v
	\label{eq: orbit, 4}
	\end{equation}
From \eqref{eq: orbit, 4} it follows that
$u\in\mathcal{O}(v,g\in G,f(g)v)$ and therefore
\[
\mathcal{O}(u,g\in G,f(g)u)\subseteq\mathcal{O}(v,g\in G,f(g)v)
\]
This completes the proof.
		\end{proof}

		 \begin{theorem}
		 \label{theorem: direct product of representations}
Suppose $f_1$ is representation of group $G$ in set $M_1$ and
$f_2$ is representation of group $G$ in set $M_2$.
Then we introduce
\AddIndex{direct product of representations
$f_1$ and $f_2$ of group}
{direct product of representations of group}
	\begin{align*}
f&=f_1\otimes f_2:G\rightarrow M_1\otimes M_2\\
f(g)&=(f_1(g),f_2(g))
	\end{align*}
		\end{theorem}
		\begin{proof}
To show that $f$ is a representation,
it is enough to prove that $f$ satisfies the definition
\ref{definition: left-side representation of group}.
\[f(e)=(f_1(e),f_2(e))=(\delta_1,\delta_2)=\delta\]
	\begin{align*}
f(ab)u&=(f_1(ab)u_1,f_2(ab)u_2)\\
&=(f_1(a)(f_1(b)u_1),f_2(a)(f_2(b)u_2))\\
&=f(a)(f_1(b)u_1,f_2(b)u_2)\\
&=f(a)(f(b)u)
	\end{align*}
		\end{proof}

			\section{Single Transitive Representation of Group}
			\label{section: Single Transitive Representation of Group}

	 \begin{definition}
We call
\AddIndex{kernel of inefficiency of representation of group}
{kernel of inefficiency of representation of group} $G$
a set \[K_f=\{g\in G:f(g)=\delta\}\]
If $K_f=\{e\}$ we call representation of group $G$
\AddIndex{effective}{effective representation of group}.
	 \qed
	 \end{definition}

		\begin{theorem}	
A kernel of inefficiency is a subgroup of group $G$.
		\end{theorem}
		\begin{proof}
Assume
$f(a_1)=\delta$ and $f(a_2)=\delta$. Then
\[f(a_1a_2)u=f(a_1)(f(a_2)u)=u\]
\[f(a^{-1})=f^{-1}(a)=\delta\]
		\end{proof}

If an action is not effective we can switch to an effective one
by changing group $G_1=G|K_f$
using factorization by the kernel of inefficiency.
This means that we can study only an effective action.

	 \begin{definition}
We call a representation of group
\AddIndex{transitive}{transitive representation of group}
if for any $a, b \in V$ exists such $g$ that
\[a=f(g)b\]
We call a representation of group
\AddIndex{single transitive}{single transitive representation of group}
if it is transitive and effective.
	 \qed
	 \end{definition}

		\begin{theorem}	
Representation is single transitive if and only if for any $a, b \in V$
exists one and only one $g\in G$ such that $a=f(g)b$
		\end{theorem}

	 \begin{definition}
We call a space $V$
\AddIndex{homogeneous space of group}{homogeneous space of group} $G$
if we have single transitive representation of group $G$ on $V$.
	 \qed
	 \end{definition}

		\begin{theorem}	%
		\label{theorem: single transitive representation of group}
If we define a single transitive representation $f$ of the group $G$ on the manifold $A$
then we can uniquely define coordinates on $A$ using coordinates on the group $G$.

If $f$ is a covariant representation than $f(a)$ is equivalent to the left shift $L(a)$ on the group $G$.
If $f$ is a contravariant representation than $f(a)$ is equivalent to the right shift $R(a)$ on the group $G$.
		\end{theorem}
		\begin{proof}
We select a point $v\in A$
and define coordinates of a point $w\in A$
as coordinates of the transformation $a$ such that $w=f(a) v$.
Coordinates defined this way are unique
up to choice of an initial point $v\in A$
because the action is effective.

If $f$ is a covariant representation we will use the notation
\[f(a)v=av\]
Because the notation
\[f(a)(f(b)v)=a(bv)=(ab)v=f(ab)v\]
is compatible with the group structure we see that the covariant representation $f$ is equivalent to the left shift.

If $f$ is a contravariant representation we will use the notation
\[f(a)v=va\]
Because the notation
\[f(a)(f(b)v)=(vb)a=v(ba)=f(ba)v\]
is compatible with the group structure we see that the contravariant representation $f$ is equivalent to the right shift.
		\end{proof}

		\begin{theorem}	%
		\label{theorem: shifts on group commuting}
Left and right shifts on group $G$ are commuting.
		\end{theorem}
		\begin{proof}
This is the consequence of the associativity on the group $G$
\[(L(a) R(b))c = a(cb)=(ac)b=(R(b) L(a))c\]
		\end{proof}

		\begin{theorem}	%
		\label{theorem: two representations of group}
If we defined a single transitive representation $f$ on the manifold $A$
then we can uniquely define a single transitive representation $h$
such that diagram
	\[
\xymatrix{
M\ar[rr]^{h(a)}\ar[d]^{f(b)} & & M\ar[d]^{f(b)}\\
M\ar[rr]_{h(a)}& &M
}
	\]
is commutative for any $a$, $b\in G$.\footnote{The theorem
\ref{theorem: two representations of group}
is really very interesting. However its
meaning becomes more clear when we apply this theorem to basis manifold,
see section \ref{section: Basis in Vector Space}.}
		\end{theorem}
		\begin{proof}
We use group coordinates for points $v\in A$.
For the simplicity we assume that $f$ is a covariant representation.
Then according to theorem \ref{theorem: single transitive representation of group}
we can write the left shift $L(a)$ instead of the transformation $f(a)$.

Let points $v_0, v\in A$. Then we can find
one and only one $a\in G$ such that
\[v=v_0 a=R(a) v_0\]
We assume
\[h(a)=R(a)\]
For some $b\in G$ we have
\[w_0=f(b)v_0=L(b) v_0\ \ \ \ w=f(b)v=L(b) v\]
According to theorem \ref{theorem: shifts on group commuting} the diagram
	\begin{equation}
\xymatrix{
v_0\ar[rr]^{h(a)=R(a)}\ar[d]^{f(b)=L(b)} & & v\ar[d]^{f(b)=L(b)}\\
w_0\ar[rr]_{h(a)=R(a)}& &w
}
	\label{Diagram: two representations of group}
	\end{equation}
is commutative.

Changing $b$ we get that $w_0$ is an arbitrary point of $A$.

We see from the diagram that if $v_0=v$ than $w_0=w$ and therefore $h(e)=\delta$.
On other hand if $v_0\neq v$ then $w_0\neq w$ because the representation $f$ is single transitive.
Therefore the representation $h$ is effective.

In the same way we can show that for given $w_0$ we can find $a$
such that $w=h(a)w_0$. Therefore the representation is single transitive.

In general the representation $f$ is not commutative and therefore
the representation $h$ is different from the representation $f$.
In the same way we can create a representation $f$ using the representation $h$.
		\end{proof}

		\begin{remark}
		\label{remark: one representation of group}
It is clear that transformations $L(a)$ and $R(a)$
are different until the group $G$ is nonabelian.
However they both are maps onto.
Theorem \ref{theorem: two representations of group} states that if both
right and left shift presentations exist on the manifold $A$
we can define two commuting representations on the manifold $A$.
The left shift or the right shift only cannot present both types of representation.
To understand why it is so let us change diagram \eqref{Diagram: two representations of group}
and assume $h(a)v_0=L(a)v_0=v$ instead of $h(a)v_0=R(a)v_0=v$ and let
us see what expression $h(a)$ has at
the point $w_0$. The diagram
\[\xymatrix{
v_0\ar[rr]^{h(a)=L(a)}\ar[d]^{f(b)=L(b)} & & v\ar[d]^{f(b)=L(b)}\\
w_0\ar[rr]_{h(a)}& &w
}\]
is equivalent to the diagram
\[\xymatrix{
v_0\ar[rr]^{h(a)=L(a)} & & v\ar[d]^{f(b)=L(b)}\\
w_0\ar[rr]_{h(a)}\ar[u]_{f^{-1}(b)=L(b^{-1})}& &w
}\]
and we have $w=bv=bav_0=bab^{-1}w_0$. Therefore
\[h(a)w_0=(bab^{-1})w_0\] We see that the representation of $h$
depends on its argument. \qed  
		\end{remark}

			\section{Linear Representation of Group}
			\label{section: Linear Representation of Group}

Suppose we introduce additional structure on set $M$.
Then we create an additional requirement for the representation of group.

Since we introduce continuity on set $M$,
we suppose that transformation
	\[
u'=f(a)u
	\]
is continuous in $u$. 
Therefore, we get
\[\left|\frac { \partial u'} {\partial u}\right|\neq 0\]

Suppose $M$ is a group. Then representations of
left and right shifts have great importance.

	 \begin{definition}
Let $M$ be vector space $\mathcal{V}$ over field $F$.
We call the representation of group $G$ in vector space $\mathcal{V}$
\AddIndex{linear representation}{linear representation of group}
if $f(a)$ is homomorphism of space $\mathcal{V}$ for any $a\in G$.
	 \qed
	 \end{definition}

		\begin{remark}
		\label{remark: row and column vectors, vector space}

Let transformation $f(a)$ be linear
homogeneous transformation.
$f_\gamma^\beta(a)$ are elements of the matrix of
transformation. We usually assume that the lower index
enumerates rows in the matrix and the upper index enumerates columns.

According to the matrix product rule we can
present coordinates of a vector as a row of a matrix.
We call such vector a \AddIndex{row vector}{row vector}.
We can also study a vector whose coordinates
form a column of a matrix.
We call such vector a \AddIndex{column vector}{column vector}.

Left-side linear representation in column vector space
\[
u'=f(a)u\ \ \ \ u'_\alpha=f_\alpha^\beta(a)u_\beta\ \ \ \ a\in G
\]
is covariant representation
\[
u''_\gamma=f_\gamma^\beta(ba) u_\beta
=f_\gamma^\alpha(b)(f_\alpha^\beta(a) u_\beta)
=(f_\gamma^\alpha(b)f_\alpha^\beta(a)) u_\beta
\]

Left-side linear representation in row vector space
\[
u'=f(a)u\ \ \ \ u'^\alpha=f^\alpha_\beta(a)u^\beta\ \ \ \ a\in G
\]
is contravariant representation
\[
u''^\gamma=f^\gamma_\beta(ba) u^\beta
=f^\gamma_\alpha(b)(f^\alpha_\beta(a) u^\beta)
=(f^\alpha_\beta(a)f^\gamma_\alpha(b)) u^\beta
\]

Right-side representation in column vector space
\[
u'=uf(a)\ \ \ \ u'_\alpha=u_\beta f_\alpha^\beta(a)\ \ \ \ a\in G
\]
is contravariant representation
\[
u''_\gamma=u_\beta f_\gamma^\beta(ab)
=(u_\beta f_\alpha^\beta(a))f_\gamma^\alpha(b)
=u_\beta(f_\gamma^\alpha(b)f_\alpha^\beta(a)) 
\]

Right-side representation in row vector space
\[
u'=uf(a)\ \ \ \ u'^\alpha=u^\beta f^\alpha_\beta(a)\ \ \ \ a\in G
\]
is covariant representation
\[
u''^\gamma=u^\beta f^\gamma_\beta(ab)
=(u^\beta f^\alpha_\beta(a))f^\gamma_\alpha(b)
=u^\beta(f^\alpha_\beta(a)f^\gamma_\alpha(b)) 
\]
		\qed
		\end{remark}

		\begin{remark}
		\label{remark: covariant and contravariant representations}
Studying linear representations
we clearly use tensor notation.
We can use only upper index and notation $u{}^\bullet{}_\alpha^.$ instead of
$u_\alpha$. Then we can write the transformation of this object in the form
\[u'{}^\bullet{}_\alpha^.=f{}^\bullet{}_\alpha^.{}{}_\bullet{}^\beta_. u{}^\bullet{}_\beta^.\]
This way we can hide the difference between covariant and
contravariant representations. This similarity goes
as far as we need.
		\qed
		\end{remark}

%% file: Basis.English.tex
\def\texBasis{}
\ifx\PrintBook\undefined
\else
				\chapter{Basis Manifold}
\fi

			\section{Basis in Vector Space}
			\label{section: Basis in Vector Space}

Assume we have vector space \symb{\mathcal{V}}1{V}
and contravariant right-side effective linear
representation of group
\symb{G(\mathcal{V})}0{GV}
$G=\ShowSymbol{GV}$.
We usually call group $G(\mathcal{V})$
\AddIndex{symmetry group}{symmetry group}.
Without loss of generality we identify element $g$ of group $G$
with corresponding transformation of representation
and write its action on vector $v\in\mathcal{V}$ as $vg$.

This point of view allows introduction of two types of coordinates
for element $g$ of group $G$. We can either use coordinates
defined on the group, or introduce coordinates as
elements of the matrix of the corresponding transformation.
The former type of coordinates is more effective when we study
properties of group $G$. The latter type of coordinates
contains redundant data;
however, it may be more convenient
when we study representation of group $G$.
The latter type of coordinates is called
\AddIndex{coordinates of representation}{coordinates of representation}.

A maximal set of linearly independent vectors
\symb{\Basis{e}}0{Basis e}
$\ShowSymbol{Basis e}=<e_{(i)}>$
is called a \AddIndex{basis}{Basis}.
In case when we want to show clearly that this is the basis
in vector $\mathcal{V}$ we use notation
\symb{\Basis{e}_\mathcal{V}}1{basis in V}.

Any homomorphism of the vector space maps one basis into another. 
Thus we can extend a covariant representation of the symmetry group
to the set of bases.
We write the action of element $g$ of group $G$ on basis $\Basis{e}$
as \symb{R(g)\Basis{e}}1{active transformation}.
However not every two bases can be mapped by a transformation
from the symmetry group
because not every nonsingular linear transformation belongs to
the representation of group $G$. Therefore, we can present
the set of bases
as a union of orbits of group $G$.

Properties of basis depend on the symmetry group.
We can select basis $\Basis{e}$ vectors of which are
in a relationship which is invariant relative to symmetry group.
In this case all bases from orbit
$\mathcal{O}(\Basis{e},g\in G,R(g)\Basis{e})$
have vectors which satisfy the same relationship. Such a basis we call
\AddIndex{\Gbasis}{G-basis}.
In each particular case we need to prove the existence of a basis
with certain properties. If such a basis does not exist
we can choose an arbitrary basis.

	 \begin{definition}
	 \label{definition: basis manifold of vector space}
We call orbit $\mathcal{O}(\Basis{e},g\in G,R(g)\Basis{e})$
of the selected basis $\Basis{e}$
the \AddIndex{basis manifold
\symb{\mathcal{B}(\mathcal{V})}1{basis manifold of vector space}
of vector space}{basis manifold of vector space}
$\mathcal{V}$.
	 \qed
	 \end{definition}

		\begin{theorem}
		\label{theorem: basis manifold of vector space}
Representation of group $G$ on basis manifold
is single transitive representation.
		\end{theorem}
		\begin{proof}
According to definition \ref{definition: basis manifold of vector space}
at least one transformation of representation is defined for any two bases.
To prove this theorem it is sufficient to show
that this transformation is unique.

Consider elements $g_1$, $g_2$ of group $G$ and a basis $\Basis{e}$ such that
	\begin{equation}
R_{g_1}\Basis{e}=R_{g_2}\Basis{e}
	\label{eq: two transformations on basis manifold, 1}
	\end{equation}
From \eqref{eq: two transformations on basis manifold, 1} it follows that
	\begin{equation}
R_{g_2^{-1}}R_{g_1}\Basis{e}=R_{g_1g_2^{-1}}\Basis{e}=\Basis{e}
	\label{eq: two transformations on basis manifold, 2}
	\end{equation}
Because any vector has a unique expansion
relative to basis $\Basis{e}$ it follows from
\eqref{eq: two transformations on basis manifold, 2}
that $R_{g_1g_2^{-1}}$ is an identical transformation
of vector space $\mathcal{V}$. $g_1=g_2$ because representation
of group $G$ is effective on vector space $\mathcal{V}$.
Statement of the theorem follows from this.
		\end{proof}

Theorem \ref{theorem: basis manifold of vector space}
means that the basis manifold $\mathcal{B}(\mathcal{V})$
is a homogenous space of group $G$.
We constructed contravariant right-side single transitive linear
representation of group $G$ on the basis manifold.
Such representation is called
\AddIndex{active representation}{active representation}.
A corresponding transformation on the basis manifold is called
\AddIndex{active transformation}{active transformation}
(\citeBib{Korn})
because the homomorphism of the vector space induced this transformation.

According to theorem \ref{theorem: single transitive representation of group}
because basis manifold $\mathcal{B}(\mathcal{V})$
is a homogenous space of group $G$
we can introduce on  $\mathcal{B}(\mathcal{V})$
two types of coordinates defined on group $G$.
In both cases coordinates of basis $\Basis{e}$ are
coordinates of the homomorphism mapping a fixed basis $\Basis{e}_0$
to the basis $\Basis{e}$.
Coordinates of representation are called
\AddIndex{standard coordinates of basis}
{standard coordinates of basis}.
We can show that standard coordinates
\symb{e^i_k}1{standard coordinates of basis}
of basis $\Basis{e}$ for certain value of $k$
are coordinates of vectors
\symb{\Vector e_k}0{vector of basis}
$\ShowSymbol{vector of basis}\in\Basis{e}$
relative to a fixed basis $\Basis{e}_0$.

Basis $\Basis{e}$ creates coordinates in $\mathcal{V}$.
In different types of space it may be done in different ways.
In affine space if node of basis is point $A$ than
point $B$ has the same coordinates as vector $\overset{\longrightarrow}{AB}$
relative basis $\Basis{e}$.
In a general case we introduce coordinates of a vector as coordinates
relative to the selected basis.
Using only bases of type $G$ means using of specific coordinates on $\mathcal{A}_n$.
To distinguish them we call this
\AddIndex{\Gcoords}{G-coordinates}.
We also call the space $\mathcal{V}$ with such coordinates
\AddIndex{\Gspace}{GSpace}.

According to theorem \ref{theorem: two representations of group} another
representation, commuting with passive, exists on the basis manifold.
As we see from remark \ref{remark: one representation of group}
transformation of this
representation is different from a passive transformation and cannot be reduced to
transformation of space $\mathcal{V}$.
To emphasize the difference this transformation is called
a \AddIndex{passive transformation}{passive transformation}
of vector space $\mathcal{V}$
and the representation is called
\AddIndex{passive representation}{passive representation}.
We write the passive transformation of basis $\Basis{e}$,
defined by element $g\in G$,
as \symb{L(g)\Basis{e}}1{passive transformation}.

			\subsection{Basis in Affine Space}
			\label{sec:AffineSpace}

We identify vectors of the affine space
\symb{\mathcal{A}_n}1{An}
with pair of points $\overset{\longrightarrow}{AB}$.
All vectors that have a common beginning $A$ create a vector space
that we call a tangent vector space $T_A\mathcal{A}_n$.

A topology that $\mathcal{A}_n$ inherits from the map $\mathcal{A}_n\rightarrow R^n$
allows us to study smooth transformations of $\mathcal{A}_n$
and their derivatives. More particularly, the derivative
of transformation $f$ maps the vector space $T_A\mathcal{A}_n$ into $T_{f(A)}\mathcal{A}_n$.
If $f$ is linear then its derivative is the same at every point.
Introducing coordinates $A^1,...,A^n$ of a point $A\in\mathcal{A}_n$ we can write
a linear transformation as
	\begin{align}
A'^i&=P^i_j A^j + R^i &\det P\ne 0
	\label{eq: AffineTransformation}
	\end{align}
Derivative of this transformation is defined by matrix $\|P^i_j\|$.
and does not depend on point $A$. Vector
$(R^1,...,R^n)$ expresses displacement in affine space.
Set of transformations \eqref{eq: AffineTransformation} is the group Lie
which we denote as \symb{GL(\mathcal{A}_n)}1{GLAn}
and call
\AddIndex{affine transformation group}{AffineTransformationGroup}.

		\begin{definition}
\AddIndex{Affine basis}{Affine Basis}
\symb{\Basis{e}=<O,\Vector e_i>}1{Affine Basis} is set of
linear independent vectors $\Vector e_i=\overset{\longrightarrow}{OA_i}=(e^1_i,...,e^n_i)$
with common start point $O=(O^1,...,O^n)$. \qed
		\end{definition}

		\begin{definition}
\AddIndex{Basis manifold \symb{\mathcal{B}(\mathcal{A}_n)}1{BAn}
of affine space}{Basis Manifold, Affine Space} 
is set of bases of this space.
		\qed
		\end{definition}

An active transformation is called
\AddIndex{affine transformation}{affine transformation}.
A passive transformation is called
\AddIndex{quasi affine transformation}{quasi affine transformation}.

If we do not concern about starting point of a vector we see little different
type of space which we call central affine space
\symb{\mathcal{CA}_n}1{CAn}.
In the central affine space we can identify all
tangent spaces and denote them $T\mathcal{CA}_n$. If we assume that the start point
of vector is origin $O$ of coordinate system in space then we can identify any
point $A\in\mathcal{CA}_n$ with the vector $a=\overset{\longrightarrow}{OA}$. This leads
to identification of $\mathcal{CA}_n$ and $T\mathcal{CA}_n$.
Now transformation is simply map
	\begin{align*}
a'^i&=P^i_j a^j & \det P\ne 0
	\end{align*}
and such transformations build up Lie group $GL_n$.

		\begin{definition}
\AddIndex{Central affine basis}{Central Affine Basis}
\symb{\Basis{e}=<\Vector e_i>}1{Central Affine Basis} is set of
linearly independent vectors $\Vector e_i=(e^1_i,...,e^n_i)$. \qed
		\end{definition}

		\begin{definition}
\AddIndex{Basis manifold \symb{\mathcal{B}(\mathcal{CA}_n)}1{BCAn}
of central affine space}{Basis Manifold, Central Affine Space}
is set of bases of this space. \qed
		\end{definition}

			\subsection{Basis in Euclid Space}
			\label{subsection: EuclideSpace}

When we introduce a metric in a central affine space we get a new geometry because
we can measure a distance and a length of vector. If a metric is
positive defined we call the space Euclid
\symb{\mathcal{E}_n}1{En}
otherwise we call the space pseudo Euclid
\symb{\mathcal{E}_{nm}}1{Enm}.

Transformations that preserve length
form Lie group $SO(n)$ for Euclid space and Lie group $SO(n,m)$ for pseudo
Euclid space where $n$ and $m$ number of positive and negative terms in
metrics.

		\begin{definition}
\AddIndex{Orthonornal basis}{Orthonornal Basis}
\symb{\Basis{e}=<\Vector e_i>}1{Orthonornal Basis} is set of
linearly independent vectors $\Vector e_i=(e^1_i,...,e^n_i)$ such that length
of each vector is $1$ and different vectors are orthogonal.
		\qed
		\end{definition}

We can prove existence of orthonormal basis
using Gram-Schmidt orthogonalization procedure.

		\begin{definition}
\AddIndex{Basis manifold \symb{\mathcal{B}(\mathcal{E}_n)}1{BEn}
of Euclid space}{Basis Manifold, Euclid Space}
is set of orthonornal bases of this space. \qed
		\end{definition}

A active transformation is called
\AddIndex{movement}{movement transformation}.
An passive transformation is called
\AddIndex{quasi movement}{quasi movement}.

			\section{Geometrical Object of Vector Space}
			\label{section: Geometrical Object of Vector Space}

An active transformation changes bases and vectors uniformly
and coordinates of vector relative basis do not change.
A passive transformation changes only the basis and it leads to change
of coordinates of vector relative to basis.

Let passive transformation $L(a)\in G$ defined by matrix
$(a^i_j)$ maps
basis $\Basis{e}=<e_i>\in\mathcal{B}(\mathcal{V})$
into basis $\Basis{e}'=<e'_i>\in\mathcal{B}(\mathcal{V})$
	\begin{equation}
	\label{eq: active transformation of vector space}
e'_j=a^i_je_i
	\end{equation}
Let vector $v\in\mathcal{V}$ have expansion
	\begin{equation}
	\label{eq: vector expansion in vector space, basis f}
v=v^ie_i
	\end{equation}
relative to basis $\Basis{e}$ and have expansion
	\begin{equation}
	\label{eq: vector expansion in vector space, basis fprim}
v=v'^ie'_i
	\end{equation}
relative to basis $\Basis{e}'$.
From \eqref{eq: active transformation of vector space}
and \eqref{eq: vector expansion in vector space, basis fprim} it follows that
	\begin{equation}
	\label{eq: vector expansion in vector space, basis f, 1}
v=v'^ja^i_je_i
	\end{equation}
Comparing \eqref{eq: vector expansion in vector space, basis f}
and \eqref{eq: vector expansion in vector space, basis f, 1} we get
	\begin{equation}
	\label{eq: coordinate transformation, 1}
v^i=v'^ja^i_j
	\end{equation}
Because $a^i_j$ is nonsingular matrix
we get from \eqref{eq: coordinate transformation, 1}
	\begin{equation}
	\label{eq: coordinate transformation}
v'^i=v^ja^{-1}{}^i_j
	\end{equation}
Coordinate transformation \eqref{eq: coordinate transformation}
does not depend on vector $v$ or basis $\Basis{e}$, but is
defined only by coordinates of vector $v$
relative basis $\Basis{e}$.

Suppose we select basis $\Basis{e}$. Then the set of coordinates
\symb{(v^i)}1{coordinates in vector space}
relative to this basis forms a vector space
\symb{\tilde{\mathcal{V}}}1{coordinate vector space}
isomorphic to vector space $\mathcal{V}$.
This vector space
is called \AddIndex{coordinate vector space}{coordinate vector space}.
This isomorphism
is called \AddIndex{coordinate isomorphism}{coordinate isomorphism}.
Denote by
\symb{\Vector\delta_k=(\delta^i_k)}1
{image of vector e_k, coordinate vector space}
the image of vector $e_k\in\Basis{e}$ under this .
		\begin{theorem}
		\label{theorem: coordinate transformations form representation, vector space}
Coordinate transformations \eqref{eq: coordinate transformation}
form right-side contravariant effective linear
representation of group $G$ which
is called \AddIndex{coordinate representation}
{coordinate representation, vector space}.
		\end{theorem}
		\begin{proof}
Suppose we have two consecutive passive transformations
$L(a)$ and $L(b)$. Coordinate transformation \eqref{eq: coordinate transformation}
corresponds to passive transformation $L(a)$.
Coordinate transformation
	\begin{equation}
	\label{eq: coordinate transformation, b}
v''^k=v'^ib^{-1}{}^k_i
	\end{equation}
corresponds to passive transformation $L_b$. Product
of coordinate transformations \eqref{eq: coordinate transformation}
and \eqref{eq: coordinate transformation, b} has form
	\begin{equation}
	\label{eq: coordinate transformation, ba}
v''^k=v^ja^{-1}{}^i_jb^{-1}{}^k_i=v^j(ba)^{-1}{}^k_j
	\end{equation}
and is coordinate transformation
corresponding to passive transformation $L_{ba}$.
It proves that coordinate transformations
form contravariant right-side linear
representation of group $G$.

Suppose coordinate transformation does not change vectors $\delta_k$.
Then unit of group $G$ corresponds to it because representation
is single transitive. Therefore,
coordinate representation is effective.
		\end{proof}

Let homomorphism of group $G$ to
the group of passive transformations
of vector space $\mathcal{W}$ be coordinated with symmetry group
of vector space $\mathcal{V}$.
This means that passive transformation $L(a)$ of vector space $\mathcal{W}$
corresponds to passive transformation $L(a)$ of vector space $\mathcal{V}$.
	\begin{equation}
	\label{eq: passive transformation of vector space W}
E'_\alpha=A^\beta_\alpha(a)E_\beta
	\end{equation}
Then coordinate transformation in $\mathcal{W}$ gets form
	\begin{equation}
	\label{eq: coordinate transformation, W}
w'^\alpha=w^\beta A(a^{-1})^\alpha_\beta=w^\beta A(a)^{-1}{}^\alpha_\beta
	\end{equation}

		\begin{definition}
		\label{definition: geometrical object, coordinate representation}
Orbit
\symb{\mathcal{O}((w,\Basis{e}_\mathcal{V}),a\in G,
(wA(a)^{-1},L(a)\Basis{e}_\mathcal{V}))}2
{geometrical object, coordinate vector space}%
is called
\AddIndex{geometrical object in coordinate representation}
{geometrical object, coordinate vector space}
defined in vector space $\mathcal{V}$.
For any basis
$\Basis{e}'_\mathcal{V}=L(a)\Basis{e}_\mathcal{V}$
corresponding point \eqref{eq: coordinate transformation, W} of orbit defines
\AddIndex{coordinates of geometrical object}
{coordinates of geometrical object, coordinate vector space}
relative basis $\Basis{e}'_\mathcal{V}$.
		\qed
		\end{definition}

		\begin{definition}
		\label{definition: geometrical object}
Orbit
\symb{\mathcal{O}((w,\Basis{e}_\mathcal{W},\Basis{e}_\mathcal{V}),
a\in G,(wA(a)^{-1},L(a)\Basis{e}_\mathcal{W},
L(a)\Basis{e}_\mathcal{V}))}2
{geometrical object, vector space}%
is called
\AddIndex{geometrical object}{geometrical object, vector space}
defined in vector space $\mathcal{V}$.
For any basis
$\Basis{e}'_\mathcal{V}=L(a)\Basis{e}_\mathcal{V}$
corresponding point \eqref{eq: coordinate transformation, W} of orbit defines
\AddIndex{coordinates of a geometrical object}{coordinates of geometrical object, vector space}
relative to basis $\Basis{e}'_\mathcal{V}$
and the corresponding vector
\[
w=w'^\alpha E'_\alpha
\]
is called
\AddIndex{representative of geometrical object}
{representative of geometrical object, vector space}
in basis $\Basis{e}'_\mathcal{V}$.
		\qed
		\end{definition}

We also say that $w$ is
a \AddIndex{geometrical object of type $A$}
{geometrical object of type A, vector space}

Since a geometrical object is an orbit of representation, we see that
according to theorem \ref{theorem: proper definition of orbit}
the definition of the geometrical object is a proper definition.

Definition \ref{definition: geometrical object, coordinate representation}
introduces a geometrical object in coordinate space.
We assume in definition \ref{definition: geometrical object}
that we selected a basis in vector space $\mathcal{W}$.
This allows using a representative of the geometrical object
instead of its coordinates.

		\begin{theorem}[\AddIndex{invariance principle}
		{invariance principle, vector space}]
		\label{theorem: invariance principle}
Representative of geometrical object does not depend on selection
of basis $\Basis{e}'_\mathcal{V}$.
		\end{theorem}
		\begin{proof}
To define representative of geometrical object,
we need to select basis $\Basis{e}_\mathcal{V}$,
basis
$\Basis{e}_\mathcal{W}=(E_\alpha)$
and coordinates of geometrical object
$w^\alpha$. Corresponding representative of geometrical object
has form
\[
w=w^\alpha E_\alpha
\]
Suppose we map basis $\Basis{e}_\mathcal{V}$
to basis $\Basis{e}'_\mathcal{V}$
by passive transformation $L(a)$.
According building this forms passive transformation
\eqref{eq: passive transformation of vector space W}
and coordinate transformation
\eqref{eq: coordinate transformation, W}. Corresponding
representative of geometrical object has form
\[
w'=w'^\alpha E'_\alpha
=w^\beta A(a)^{-1}{}^\alpha_\beta A^\gamma_\alpha(a)E_\gamma
=w'^\beta E'_\beta=w
\]
Therefore representative of geometrical object
is invariant relative selection of basis.
		\end{proof}

		\begin{definition}
		\label{definition: sum of geometrical objects}
Let
\[w_1=w_1^\alpha E_\alpha\]
\[w_2=w_2^\alpha E_\alpha\]
be geometrical objects of the same type
defined in vector space $\mathcal{V}$.
Geometrical object
\[w=(w_1^\alpha+w_2^\alpha)E_\alpha\]
is called \AddIndex{sum
\[w=w_1+w_2\]
of geometrical objects}{sum of geometrical objects, vector space}
$w_1$ and $w_2$.
		\qed
		\end{definition}

		\begin{definition}
		\label{definition: product of geometrical object and constant}
Let
\[w_1=w_1^\alpha E_\alpha\]
be geometrical object
defined in vector space $\mathcal{V}$ over field $F$.
Geometrical object
\[w_2=(kw_1^\alpha) E_\alpha\]
is called \AddIndex{product
\[w_2=kw_1\]
of geometrical object $w_1$ and constant $k\in F$}
{product of geometrical object and constant, vector space}.
		\qed
		\end{definition}

		\begin{theorem}
Geometrical objects of type $A$
defined in vector space $\mathcal{V}$ over field $F$
form vector space over field $F$.
		\end{theorem}
		\begin{proof}
The statement of the theorem follows from immediate verification
of the properties of vector space.
		\end{proof}

%% file: Biblio.English.tex
\OpenBiblio

\BiblioItem{texIntro}{Einstein: Geometry and Experience}
{
Einstein, Geometry and Experience, (1921)
}%

\BiblioItem{texGenRelativity}{Ghez}
{
Ghez et al.,
The First Measurement of Spectral Lines in a Short-Period Star Bound to the Galaxy's Central Black Hole: A Paradox of Youth,
\href{http://www.journals.uchicago.edu/ApJ/journal/issues/ApJL/v586n2/16990/brief/16990.abstract.html}{ApJL, 586, L127} (2003),
eprint \href{http://arxiv.org/abs/astro-ph/0302299}{arXiv:astro-ph/0302299} (2003)
}%

\BiblioItem{texGenRelativity}{Schodel}
{
R. Sch\"odel et al.,
A star in a 15.2-year orbit around the supermassive black hole at the centre of the Milky Way,
\href{http://www.nature.com/cgi-taf/DynaPage.taf?file=/nature/journal/v419/n6908/abs/nature01121_fs.html}{Nature 419, 694} (2002)
}%

\BiblioItem{texAffine}{Mielke}
{
Eckehard W. Mielke, Affine generalization of the Komar complex of general relativity,
\href{http://prola.aps.org/searchabstract/PRD/v63/i4/e044018}{Phys. Rev. D 63, 044018} (2001)
}%

\BiblioItem{texAffine}{Obukhov}
{
Yu. N. Obukhov and J. G. Pereira, Metric-affine approach to teleparallel gravity,
\href{http://scitation.aip.org/getabs/servlet/GetabsServlet?prog=normal&id=PRVDAQ000067000004044016000001&idtype=cvips&gifs=Yes}
{Phys. Rev. D 67, 044016} (2003),
eprint \href{http://arxiv.org/abs/gr-qc/0212080}{arXiv:gr-qc/0212080} (2002)
}%

\BiblioItem{texAffine}{Sardanashvily}
{
Giovanni Giachetta, Gennadi Sardanashvily, Dirac Equation in Gauge and Affine-Metric Gravitation Theories,
eprint \href{http://arxiv.org/abs/gr-qc/9511035}{arXiv:gr-qc/9511035} (1995)
}%

\BiblioItem{texAffine}{Gauge}
{
Frank Gronwald and Friedrich W. Hehl, On the Gauge Aspects of Gravity, eprint
\href{http://arxiv.org/abs/gr-qc/9602013}{arXiv:gr-qc/9602013} (1996)
}%

\BiblioItem{texAffine}{Neeman}
{
Yuval Neeman, Friedrich W. Hehl, Test Matter in a Spacetime with Nonmetricity, eprint
\href{http://arxiv.org/abs/gr-qc/9604047}{arXiv:gr-qc/9604047} (1996)
}%

\BiblioItem{texAffine}{0405.027}
{
Aleks Kleyn,
Reference Frame in General Relativity,
eprint \href{http://arxiv.org/abs/gr-qc/0405027}{arXiv:gr-qc/0405027} (2004)
}%

\BiblioItem{texTidal,texAffine}{torsion}
{
F. W. Hehl, P. von der Heyde, G. D. Kerlick, and J. M. Nester,
General relativity with spin and torsion: Foundations and prospects,
\href{http://prola.aps.org/abstract/RMP/v48/i3/p393_1}{Rev. Mod. Phys. 48, 393} (1976)
}%

\BiblioItem{texTidal,texNewton}{Megged}
{
O. Megged, Post-Riemannian Merger of Yang-Mills Interactions with Gravity,
eprint \href{http://arxiv.org/abs/hep-th/0008135}{arXiv:hep-th/0008135} (2001)
}%


\BiblioItem{texNewton}{gr-qc-9604027}
{
Yu.N. Obukhov, E.J. Vlachynsky, W. Esser, R. Tresguerres and F.W. Hehl,
An exact solution of the metric-affine gauge theory with dilation, shear, and spin charges,
eprint \href{http://arxiv.org/abs/gr-qc/9604027}{arXiv:gr-qc/9604027} (1996)
}%

\BiblioItem{texLagrange}{Weinberg}
{
Steven Weinberg. The Quantum Theory of Fields. Cambridge university press.
}%

\BiblioItem{texLagrange}{Reinhardt}
{
Greiner Reinhardt. Field Quantization. Springer.
}%

\BiblioItem{texLagrange}{Landau}
{
L. D. Landau, E. M. Lifshich, The classical theory of fields.
Oxford, New York, Pergamon Press
}%

\BiblioItem{texTidal}{Wheeler}
{
Ignazio Ciufolini, John Wheeler. Gravitation and Inertia.
Princeton university press.
}%

\BiblioItem{texTidal}{0405.028}
{
Aleks Kleyn, Metric-Affine Manifold,
eprint \href{http://arxiv.org/abs/gr-qc/0405028}{arXiv:gr-qc/0405028} (2004)
}%

\BiblioItem{texTidal}{Anderson02}
{
J. D. Anderson, P. A. Laing, E. L. Lau, A. S. Liu, M. M. Nieto, and S. G. Turyshev,
Study of the anomalous acceleration of Pioneer 10 and 11,
\href{http://prola.aps.org/searchabstract/PRD/v65/i8/e082004}{Phys. Rev. D 65, 082004, 50 pp.}, (2002),
eprint \href{http://arxiv.org/abs/gr-qc/0104064}{arXiv:gr-qc/0104064} (2001)
}%

\BiblioItem{texTidal}{Anderson98}
{
J. D. Anderson, P. A. Laing, E. L. Lau, A. S. Liu, M. M. Nieto, and S. G. Turyshev,
Indication, from Pioneer 10/11, Galileo, and Ulysses Data, of an Apparent Anomalous, Weak, Long-Range Acceleration,
\href{http://prola.aps.org/abstract/PRL/v81/i14/p2858_1}{Phys. Rev. Lett. 81, 2858}, (1998),
eprint \href{http://arxiv.org/abs/gr-qc/9808081}{arXiv:gr-qc/9808081} (1998)
}%


\BiblioItem{texReferenceFrame,texAlgebraBundle}{Serge Lang}
{
Serge Lang,
Algebra, Springer, 2002
}%

\BiblioItem{texAlgebraBundle,texDrcRepresentation}{Burris Sankappanavar}
{
S. Burris, H.P. Sankappanavar,
A Course in Universal Algebra, Springer-Verlag (March, 1982),
\\eprint
\href{http://www.math.uwaterloo.ca/~snburris/htdocs/ualg.html}
{http://www.math.uwaterloo.ca/~snburris/htdocs/ualg.html}
\\(The Millennium Edition)
}%


\BiblioItem{texAffine,texRepresentation,texBasis,texDrcBasis,texLinearMap,texVectorSpace}{Rashevsky}
{
P. K. Rashevsky, Riemann Geometry and Tensor Calculus,\\
Moscow, Nauka, 1967
}%

\BiblioItem{texDrcBasis,texBasis}{Korn}
{
Granino A. Korn, Theresa M. Korn,
Mathematical Handbook for Scientists and Engineer,
McGraw-Hill Book Company, New York, San Francisco,
Toronto, London, Sydney, 1968
}%


\BiblioItem{texGenRelativity}{Tartaglia}
{
Angelo Tartaglia and Matteo Luca Ruggiero,
Angular Momentum Effects in Michelson\Hyph Morley Type Experiments,
Gen.Rel.Grav. 34, 1371-1382 (2002),\\
eprint \href{http://arxiv.org/abs/gr-qc/0110015}{arXiv:gr-qc/0110015} (2001)
}%

\BiblioItem{texGenRelativity}{Tomozawa}
{
Yukio Tomozawa, Speed of Light in Gravitational Fields, eprint
\href{http://arxiv.org/abs/astro-ph/0303047}{arXiv:astro-ph/0303047} (2004)
}%

\BiblioItem{texGenRelativity}{Magueijo}
{
Joao Magueijo,
Covariant and locally Lorentz-invariant varying speed of light theories,
\href{http://prola.aps.org/abstract/PRD/v62/i10/e103521}{Phys. Rev. D 62, 103521} (2000),
eprint \href{http://arxiv.org/abs/gr-qc/0007036}{arXiv:gr-qc/0007036} (2000)
}%

\BiblioItem{texGenRelativity}{Bassett}
{
Bruce A. Bassett, Stefano Liberati, Carmen Molina-Paris, and Matt Visser,
Geometrodynamics of variable-speed-of-light cosmologies,
\href{http://prola.aps.org/abstract/PRD/v62/i10/e103518}{Phys. Rev. D 62}, 103518 (2000),
eprint \href{http://arxiv.org/abs/astro-ph/0001441}{arXiv:astro-ph/0001441} (2000)
}%

\BiblioItem{texGenRelativity}{Straumann}
{
Lochlainn O'Raifeartaigh and Norbert Straumann,
Gauge theory: Historical origins and some modern developments,
\href{http://prola.aps.org/abstract/RMP/v72/i1/p1_1}{Rev. Mod. Phys. 72, 1} (2000)
}%

\BiblioItem{texGenRelativity}{Lammerzahl}
{
Claus L\"ammerzahl, Mark P. Haugan,
On the interpretation of Michelson\Hyph Morley experiments,
{Phys. Lett. A282 223-229} (2001),\\
eprint \href{http://arxiv.org/abs/gr-qc/0103052}{arXiv:gr-qc/0103052} (2001)
}%

\BiblioItem{texGenRelativity}{Muller}
{
Holger Muller et al.,
Modern Michelson-Morley Experiment using Cryogenic Optical Resonators,
\href{http://prola.aps.org/searchabstract/PRL/v91/i2/e020401}{Phys. Rev. Lett. 91, 020401} (2003),
eprint \href{http://arxiv.org/abs/physics/0305117}{arXiv:physics/0305117} (2000)
}%

\BiblioItem{texGenRelativity,texTidal}{Ranada}
{
Antonio F. Ranada,
Pioneer acceleration and variation of light speed: experimental situation,
eprint \href{http://arxiv.org/abs/gr-qc/0402120}{arXiv:gr-qc/0402120} (2004)
}%

\BiblioItem{texBiring,texVectorSpace}{math.QA-0208146}
{
I. Gelfand, S. Gelfand, V. Retakh, R. Wilson,
Quasideterminants,\\
eprint \href{http://arxiv.org/abs/math.QA/0208146}{arXiv:math.QA/0208146} (2002)
}%

\BiblioItem{texBiring,texVectorSpace}{q-alg-9705026}
{
I.Gelfand, V.Retakh,
Quasideterminants, I,\\
eprint \href{http://arxiv.org/abs/q-alg/9705026}{arXiv:q-alg/9705026} (1997)
}%

\BiblioItem{texVectorSpace}{Gelfand Retakh 1991}
{
I. Gelfand and V. Retakh, Determinants of Matrices over Noncommutative Rings, Funct.
Anal. Appl. 25 (1991), no. 2, 91-102
}%

\BiblioItem{texVectorSpace}{Gelfand Retakh 1992}
{
I. Gelfand and V. Retakh, A Theory of Noncommutative Determinants and Characteristic
Functions of Graphs, Funct. Anal. Appl. 26 (1992), no. 4, 1-20
}%

\BiblioItem{texVectorSpace}{hep-th-9407124}
{
I. M. Gelfand, D. Krob, A. Lascoux, B. Leclerc, V.S. Retakh and J.-Y. Thibon,
Noncommutative symmetric functions,\\
eprint \href{http://arxiv.org/abs/hep-th/9407124}{arXiv:hep-th/9407124} (1994)
}%

\BiblioItem{texVectorSpace}{Carl Faith 1}
{
Carl Faith, Algebra: Rings, Modules and Categories I,
Springer - Verlag, Berlin - Heidelberg - New York, 1973
}%



\BiblioItem{texReferenceFrame}{math.DG-0412391}
{
Aleks Kleyn,
Basis Manifold,
eprint \href{http://arxiv.org/abs/math.DG/0412391}{arXiv:math.DG/0412391} (2004)
}%

\BiblioItem{texAlgebraBundle,texBundleRelation}{0701.238}
{
Aleks Kleyn,
Lectures on Linear Algebra over Skew Field,\\
eprint \href{http://arxiv.org/abs/math.GM/0701238}{arXiv:math.GM/0701238} (2007)
}%

\BiblioItem{texBundleRelation,texPrefaceRelation}{0702.561}
{
Aleks Kleyn,
Algebra Bundle,\\
eprint \href{http://arxiv.org/abs/math.DG/0702561}{arXiv:math.DG/0702561} (2007)
}%


\BiblioItem{texPolymodule}{math.RA-0501237v1}
{
Aleks Kleyn,
Module Over Skew-Field, version 1,\\
eprint \href{http://arxiv.org/abs/math/0501237v1}{arXiv:math.RA/0501237v1} (2005)
}%

\ifx\PrintBook\Defined
\else
\BiblioItem{texVectorSpace,texAlgebraBundle}{0612.111}
{
Aleks Kleyn,
Biring of Matrices,\\
eprint \href{http://arxiv.org/abs/math.OA/0612111}{arXiv:math.OA/0612111} (2006)
}%
\fi

\BiblioItem{texHomotopy}{q-alg-9705009}
{
John C. Baez,
An Introduction to n-Categories,\\
eprint \href{http://arxiv.org/abs/q-alg/9705009}{arXiv:q-alg/9705009} (1997)
}%

\BiblioItem{texIntro}{Einstein: Isaak Newton}
{
Albert Einstein,
Isaak Newton,
Manchester Guardian, 16, 234 - 235 (1927)
}%

\BiblioItem{texPrefaceRelation}{Tolstoi about Anna Karenina}
{
Tolstoi about Anna Karenina,
in book A Karenina Companion, by C. J. G. Turner,
published by Wilfrid Laurier University Press (August 1993)
}%

\BiblioItem{texBundleRelation,texPrefaceRelation}{Cohn: Universal Algebra}
{
Paul M. Cohn,
Universal Algebra,
Springer, 1981
}%

\BiblioItem{texAlgebraBundle}{Pommaret: Partial Differential Equations}
{
J.-F. Pommaret,
Partial Differential Equations and Group Theory,
Springer, 1994
}%

\BiblioItem{texBundleRelation}{Bourbaki: Set Theory}
{
N. Bourbaki,
Theory of sets,
Springer, 2004
}%

\BiblioItem{texCalculus}{Pontryagin: Topological Group}
{
L. S. Pontryagin,
Selected Works, Volume Two, Topological Groups,
Gordon and Breach Science Publishers, 1986
}

\BiblioItem{texCalculus}{Bourbaki: Topological Vector Space}
{
N. Bourbaki,
Topological Vector Spaces, Chapters 1 - 5,
Springer, 1987
}

\BiblioItem{texCalculus}{Bourbaki: General Topology: Chapter 5 - 10}
{
N. Bourbaki,
General Topology, Chapters 5 - 10,
Springer, 1989
}

\BiblioItem{texAlgebraBundle}{Hatcher: Algebraic Topology}
{
Allen Hatcher,
Algebraic Topology,
Cambridge University Press, 2002
}

\CloseBiblio

%% file: Index.English.tex
\OpenIndex
\SetIndexSpace%
\Index{texLinearMap}
   {$1$-\drc form}%
   {1-drc form, vector spaces}%
\SetIndexSpace%
\Index{texPolymodule}
   {$(2)$\hyph vector space}%
   {(2)-vector space}%
\Index{texBundleRelation}
   {$2$\Hyph ary fibered relation}%
   {2 ary fibered relation}%
\SetIndexSpace%
\Index{texCalculus}
   {$A$\Hyph valued function}%
   {A valued function}%
\Index{texBiring}
   {$(^{\gi a}_{\gi b})$\hyph \CR quasideterminant}%
   {a b cr-quasideterminant}%
\Index{texBiring}
   {$(^{\gi a}_{\gi b})$\hyph \RC quasideterminant}%
   {a b RC-quasideterminant}%
\Index{texCalculus}
   {absolute value on skew field}%
   {absolute value on skew field}%
\Index{texBasis}
   {active representation}%
   {active representation}%
\Index{texDrcBasis}
   {active \sT representation}%
   {active representation, vector space}%
\Index{texBasis}
   {active transformation on basis manifold}%
   {active transformation}%
\Index{texBasis}
   {affine basis}%
   {Affine Basis}%
\Index{texBasis}
   {affine transformation group}%
   {AffineTransformationGroup}%
\Index{texTypeBasis}
   {affine transformation group}%
   {AffineTransformationGroup}%
\Index{texBasis}
   {affine transformation on basis manifold}%
   {affine transformation}%
\Index{texAlgebraBundle}
   {algebra bundle}%
   {algebra bundle}%
\Index{texBiring}
   {alternative representation of matrix}%
   {Alternative representation}%
\Index{texReferenceFrame}
   {anholonomic coordinate}%
   {anholonomic coordinate}%
\Index{texReferenceFrame}
   {anholonomic coordinates of vector}%
   {vector anholonomic coordinates}%
\Index{texReferenceFrame}
   {anholonomic coordinates on manifold}%
   {anholonomic coordinates on manifold}%
\Index{texReferenceFrame}
   {anholonomity object}%
   {anholonomity object}%
\Index{texAlgebraBundle}
   {antihomomorphism of group bundles}%
   {antihomomorphism of group bundles}%
\Index{texBundleRelation}
   {antisymmetric $2$\Hyph ary fibered relation}%
   {antisymmetric 2 ary fibered relation}%
\Index{texAlgebraBundle}
   {arity of operation}%
   {arity of operation}%
\Index{texLinearMap,texDrcRepresentation,texVectorSpace}
   {associative law}%
   {}%
\SubIndex{texDrcRepresentation}
   {for covariant \Ts representation}%
   {associative law for Tstar covariant representation}%
\SubIndex{texVectorSpace}
   {for \drc linear maps of vector spaces}%
   {associative law for drc linear maps of vector spaces}%
\SubIndex{texVectorSpace}
   {for \Ts vector space}%
   {associative law, Tstar vector space}%
\SubIndex{texLinearMap}
   {for twin representations}%
   {associative law for twin representations}%
\Index{texBundleRelation}
   {associative law of composition of fibered correspondences}%
   {associative law, composition of fibered correspondences}%
\Index{texAffine}
   {auto parallel line}%
   {auto parallel line}%
\SetIndexSpace%
\Index{texBundleRelation}
   {base of fibered correspondence}%
   {base of fibered correspondence}%
\Index{texAlgebraBundle}
   {base of map}%
   {base of map}%
\Index{texTypeBasis}
   {basis}%
   {}%
\SubIndex{texTypeBasis}
   {affine}%
   {Affine Basis}%
\SubIndex{texTypeBasis}
   {central affine}%
   {Central Affine Basis}%
\SubIndex{texTypeBasis}
   {orthonornal}%
   {Orthonornal Basis}%
\Index{texDrcBasis}
   {basis manifold}%
   {}%
\SubIndex{texTypeBasis}
   {of affine space}%
   {Basis Manifold, Affine Space}%
\SubIndex{texTypeBasis}
   {of central affine space}%
   {Basis Manifold, Central Affine Space}%
\SubIndex{texDrcBasis}
   {of \drc vector space}%
   {basis manifold of vector space}%
\SubIndex{texTypeBasis}
   {of Euclid space}%
   {Basis Manifold, Euclid Space}%
\Index{texBasis}
   {basis manifold of affine space}%
   {Basis Manifold, Affine Space}%
\Index{texBasis}
   {basis manifold of central affine space}%
   {Basis Manifold, Central Affine Space}%
\Index{texBasis}
   {basis manifold of Euclid space}%
   {Basis Manifold, Euclid Space}%
\Index{texBasis}
   {basis manifold of vector space}%
   {basis manifold of vector space}%
\Index{texBasis}
   {basis of vector space}%
   {Basis}%
\Index{texLieRepresentation}
   {basis vector}%
   {}%
\SubIndex{texLieRepresentation}
   {of \sT representation}%
   {basis vector of starT representation}%
\SubIndex{texLieRepresentation}
   {of \Ts representation}%
   {basis vector of Tstar representation}%
\Index{texBiring}
   {biring}%
   {biring}%
\SetIndexSpace%
\Index{texVectorSpace}
   {\subs rows \dcr vector space}%
   {subs rows dcr vector space}%
\Index{texBiring}
   {\subs row of matrix}%
   {c row}%
\Index{texBiring}
   {$c$\hyph row of matrix}%
   {c-row}%
\Index{texAffine}
   {Cartan connection}%
   {Cartan connection}%
\Index{texAffine}
   {Cartan curvature}%
   {Cartan curvature}%
\Index{texAffine}
   {Cartan derivative}%
   {Cartan derivative}%
\Index{texAffine}
   {Cartan symbol}%
   {Cartan symbol}%
\Index{texAffine}
   {Cartan transport}%
   {Cartan transport}%
\Index{texAlgebraBundle}
   {Cartesian power of bundle}%
   {Cartesian power of bundle}%
\Index{texAlgebraBundle}
   {direct product of bundles}%
   {Cartesian product of bundles}%
\Index{texAlgebraBundle}
   {direct product of total spaces}%
   {Cartesian product of total spaces}%
\Index{texHomotopy}
   {Category of \drc vector spaces}%
   {Category of Drc vector spaces}%
\Index{texBundleRelation}
   {category of fibered correspondences over diagonal}%
   {category of fibered correspondences over diagonal}%
\Index{texBundleRelation}
   {category of reduced fibered correspondences}%
   {category of reduced fibered correspondences}%
\Index{texBasis}
   {central affine basis}%
   {Central Affine Basis}%
\Index{texRepresentation}
   {column vector}%
   {column vector}%
\Index{texBundleRelation}
   {commutative diagram of correspondences}%
   {commutative diagram of correspondences}%
\Index{texDiffEq}
   {completely integrable system}%
   {completely integrable system}%
\Index{texBundleRelation}
   {composition of fibered correspondences}%
   {composition of fibered correspondences}%
\SubIndex{texVectorSpace}
   {of \drc linear equations}%
   {extended matrix, system of drc linear equations}%
\SubIndex{texVectorSpace}
   {of \rcd linear equations}%
   {extended matrix, system of rcd linear equations}%
\Index{texBundleRelation}
   {composition of reduced fibered correspondences}%
   {composition of reduced fibered correspondences}%
\Index{texBiring}
   {condition of reducibility of products}%
   {condition of reducibility of products}%
\Index{texBasis}
   {coordinate isomorphism}%
   {coordinate isomorphism}%
\Index{texVectorSpace}
   {coordinate isomorphism}%
   {coordinate isomorphism}%
\Index{texVectorSpace}
   {coordinate matrix}%
   {}%
\SubIndex{texVectorSpace}
   {of set of vectors in \dcr rows vector space}%
   {coordinate matrix of set of vectors, dcr vector space}%
\SubIndex{texVectorSpace}
   {of set of vectors in \drc rows vector space}%
   {coordinate matrix of set of vectors, drc vector space}%
\SubIndex{texVectorSpace}
   {of vector in \drc basis}%
   {coordinate matrix of vector in drc basis}%
\Index{texReferenceFrame}
   {coordinate reference frame}%
   {coordinate reference frame}%
\Index{texGenRelativity}
   {coordinate reference frame in event space}%
   {coordinate reference frame in event space}%
\Index{texDrcBasis}
   {coordinate representation in \drc vector space}%
   {coordinate representation, vector space}%
\Index{texBasis}
   {coordinate representation of group in vector space}%
   {coordinate representation, vector space}%
\Index{texBasis}
   {coordinate vector space}%
   {coordinate vector space}%
\Index{texDrcBasis}
   {coordinates of geometrical object}%
   {}%
\SubIndex{texDrcBasis}
   {in coordinate vector space}%
   {coordinates of geometrical object, coordinate vector space}%
\SubIndex{texDrcBasis}
   {in vector space}%
   {coordinates of geometrical object, vector space}%
\Index{texBasis}
   {coordinates of geometrical object}%
   {coordinates of geometrical object, vector space}%
\Index{texBasis}
   {coordinates of geometrical object in coordinate representation}%
   {coordinates of geometrical object, coordinate vector space}%
\Index{texDrcBasis}
   {coordinates of representation}%
   {coordinates of representation}%
\Index{texBasis}
   {coordinates of representation}%
   {coordinates of representation}%
\Index{texVectorSpace}
   {coordinates of set of vectors in \dcr vector space}%
   {coordinates of set of vectors, dcr vector space}%
\Index{texVectorSpace}
   {coordinates of set of vectors in \drc vector space}%
   {coordinates of set of vectors, drc vector space}%
\Index{texVectorSpace}
   {coordinates of vector in \drc basis}%
   {coordinates of vector in drc basis}%
\Index{texBiring}
   {\CR inverse element of biring}%
   {cr-inverse element}%
\Index{texVectorSpace}
   {\CR matrix group}%
   {cr-matrix group}%
\Index{texBiring}
   {\CR power}%
   {cr power}%
\Index{texBiring}
   {\CR product of matrices}%
   {cr-product of matrices}%
\Index{texVectorSpace}
   {\crd vector space}%
   {crd vector space}%
\SetIndexSpace%
\Index{texCalculus}
   {$D$\Hyph valued variable}%
   {D valued variable}%
\Index{texVectorSpace}
   {\dcr basis of \subs rows vector space}%
   {dcr basis, c rows vector space}%
\Index{texVectorSpace}
   {\dcr vector}%
   {dcr vector}%
\Index{texVectorSpace}
   {\dcr vector space}%
   {dcr vector space}%
\Index{texBiring}
   {determinant of matrix}%
   {determinant}%
\Index{texTidal}
   {deviation of trajectories}%
   {deviation of trajectories}%
\Index{texBundleRelation}
   {diagonal in bundle}%
   {diagonal in bundle}%
\Index{texBundleRelation}
   {diagram of correspondences}%
   {diagram of correspondences}%
\Index{texCalculus}
   {differentiable functions of \drc vector space to skew field $D$}%
   {differentiable functions, drc vector space to skew field}%
\Index{texCalculus}
   {differential of mapping of normed \drc vector space to valued skew field}%
   {differential, drc vector space to skew field}%
\Index{texVectorSpace}
   {dimension of \drc vector space}%
   {dimension of vector space}%
\Index{texRepresentation}
   {direct product of representations of group}%
   {direct product of representations of group}%
\Index{texAlgebraBundle}
   {direct product of representations of group bundle}%
   {direct product of representations of group bundle}%
\Index{texDrcRepresentation}
   {direct product of \Ts representations of group}%
   {direct product of representations of group}%
\Index{texLieRepresentation}
   {direct sum of representations}%
   {direct sum of representations}%
\Index{texVectorSpace}
   {distributive law}%
   {}%
\SubIndex{texVectorSpace}
   {\Ts vector space}%
   {distributive law, Tstar vector space}%
\Index{texVectorSpace}
   {\drc automorphism of vector space}%
   {automorphism of vector space}%
\Index{texVectorSpace}
   {\drc basis}%
   {}%
\SubIndex{texVectorSpace}
   {for \sups rows vector space}%
   {drc basis, r rows vector space}%
\SubIndex{texVectorSpace}
   {for vector space}%
   {drc basis, vector space}%
\Index{texVectorSpace}
   {drc coordinate vector space}%
   {drc coordinate vector space}%
\Index{texVectorSpace}
   {\drc isomorphism of vector spaces}%
   {isomorphism of vector spaces}%
\Index{texVectorSpace}
   {\drc linear map of vector spaces}%
   {drc linear map of vector spaces}%
\Index{texVectorSpace}
   {\drc linear span in vector space}%
   {linear span, vector space}%
\Index{texVectorSpace}
   {\drc linearly dependent vectors}%
   {linearly dependent, vector space}%
\Index{texReferenceFrame}
   {\drc linearly independent vector fields}%
   {linearly independent vector fields}%
\Index{texVectorSpace}
   {\drc linearly independent vectors}%
   {linearly independent, vector space}%
\Index{texVectorSpace}
   {\drc system of linear equations}%
   {system of linear equations}%
\Index{texVectorSpace}
   {\drc vector}%
   {drc vector}%
\Index{texCalculus}
   {\drc vector function}%
   {drc vector function}%
\Index{texVectorSpace}
   {\drc vector space}%
   {drc vector space}%
\Index{texVectorSpace}
   {$D\star$\hyph  vector space}%
   {Dstar vector space}%
\Index{texVectorSpace}
   {$D\star$\hyph product of vector over scalar}%
   {Dstar product of vector over scalar, vector space}%
\Index{texBiring}
   {duality principle for biring}%
   {duality principle for biring}%
\Index{texBiring}
   {duality principle for biring of matrices}%
   {duality principle for biring of matrices}%
\SetIndexSpace%
\Index{texRepresentation}
   {effective representation of group}%
   {effective representation of group}%
\Index{texVectorSpace}
   {effective representation of skew field}%
   {effective representation of skew field}%
\Index{texELie}
   {enhanced Lie group}%
   {enhanced Lie group}%
\Index{texDiffEq}
   {essential parameters}%
   {essential parameters}%
\Index{texBundleRelation}
   {extension of correspondence}%
   {extension of correspondence}%
\Index{texAffine}
   {extreme line}%
   {extreme line}%
\SetIndexSpace%
\Index{texBundleRelation}
   {fibered correspondence from $\mathcal{A}$ to $\mathcal{B}$}%
   {fibered correspondence from A to B}%
\Index{texBundleRelation}
   {fibered correspondence in $\mathcal{A}$}%
   {fibered correspondence in A}%
\Index{texBundleRelation}
   {fibered equivalence}%
   {fibered equivalence}%
\Index{texBundleRelation}
   {fibered ordering}%
   {fibered ordering}%
\Index{texBundleRelation}
   {fibered preordering}%
   {fibered preordering}%
\Index{texBundleRelation}
   {fibered subset}%
   {fibered subset}%
\Index{texNewton}
   {field-strength tensor}%
   {field-strength tensor}%
\Index{texNewton}
   {first Newton law}%
   {First Newton law}%
\Index{texAffine}
   {Frenet transport}%
   {Frenet transport}%
\Index{texCalculus}
   {function continuous with respect to set of arguments}%
   {function continuous with respect to set of arguments}%
\Index{texCalculus}
   {function of $\gi n$ $D$\Hyph valued variables}%
   {function of n D valued variables}%
\SetIndexSpace%
\Index{texTypeBasis}
   {\Gbasis}%
   {G-basis}%
\Index{texBasis}
   {\Gbasis\ of vector space}%
   {G-basis}%
\Index{texBasis}
   {\Gcoords\ of basis}%
   {G-coordinates}%
\Index{texTypeBasis}
   {\Gcoords}%
   {G-coordinates}%
\Index{texTypeBasis}
   {\Gspace}%
   {GSpace}%
\Index{texDrcBasis}
   {geometrical object}%
   {}%
\SubIndex{texDrcBasis}
   {defined in vector space}%
   {geometrical object, vector space}%
\SubIndex{texDrcBasis}
   {in coordinate representation defined in vector space}%
   {geometrical object, coordinate vector space}%
\SubIndex{texDrcBasis}
   {of type $A$}%
   {geometrical object of type A, vector space}%
\Index{texBasis}
   {geometrical object in coordinate representation}%
   {geometrical object, coordinate vector space}%
\Index{texBasis}
   {geometrical object in vector space}%
   {geometrical object, vector space}%
\Index{texBasis}
   {geometrical object of type $A$ in vector space}%
   {geometrical object of type A, vector space}%
\Index{texGroupRing}
   {group algebra}%
   {group algebra}%
\Index{texAlgebraBundle}
   {group bundle}%
   {group bundle}%
\Index{texBasis}
   {\Gspace}%
   {GSpace}%
\SetIndexSpace%
\Index{texBiring}
   {Hadamard inverse of matrix}%
   {Hadamard inverse of matrix}%
\Index{texReferenceFrame}
   {holonomic coordinates of vector}%
   {vector holonomic coordinates}%
\Index{texAlgebraBundle}
   {homogeneous bundle of group bundle}%
   {homogeneous bundle of group bundle}%
\Index{texDrcRepresentation}
   {homogeneous space of group}%
   {homogeneous space of group}%
\Index{texRepresentation}
   {homogeneous space of group}%
   {homogeneous space of group}%
\Index{texAlgebraBundle}
   {homomorphism of algebra bundles}%
   {homomorphism of algebra bundles}%
\Index{texAlgebraBundle}
   {homomorphism of group bundles}%
   {homomorphism of group bundles}%
\SetIndexSpace%
\Index{texLieRepresentation}
   {infinitesimal generator}%
   {infinitesimal generator}%
\Index{texLinearLie}
   {infinitesimal generators of group Lie}%
   {infinitesimal generators of group Lie}%
\Index{texDrcBasis}
   {invariance principle}%
   {invariance principle}%
\Index{texBasis}
   {invariance principle in vector space}%
   {invariance principle, vector space}%
\Index{texBundleRelation}
   {inverse fibered correspondence}%
   {inverse fibered correspondence}%
\Index{texBundleRelation}
   {inverse reduced fibered correspondence}%
   {inverse reduced fibered correspondence}%
\Index{texAlgebraBundle}
   {isomorphism of algebra bundles}%
   {isomorphism of algebra bundles}%
\SetIndexSpace%
\Index{texRepresentation}
   {kernel of inefficiency of representation of group}%
   {kernel of inefficiency of representation of group}%
\Index{texAlgebraBundle}
   {kernel of inefficiency of representation of group bundle}%
   {kernel of inefficiency of representation of group bundle}%
\Index{texDrcRepresentation}
   {kernel of inefficiency of \Ts representation of group $G$}%
   {kernel of inefficiency of representation of group}%
\Index{texDiffProperty}
   {Killing equation}%
   {Killing equation}%
\Index{texDiffProperty}
   {Killing equation of second type}%
   {Killing equation second type}%
\Index{texDiffProperty}
   {Killing vector of second type}%
   {Killing vector second type}%
\Index{texBiring}
   {Kronecker symbol}%
   {Kronecker symbol}%
\SetIndexSpace%
\Index{texLie}
   {left invariant vector field}%
   {left invariant vector}%
\Index{texVectorSpace}
   {left module}%
   {left module}%
\Index{texDrcRepresentation}
   {left shift}%
   {left shift}%
\Index{texRepresentation}
   {left shift on group}%
   {left shift, group}%
\Index{texAlgebraBundle}
   {left shift on group bundle}%
   {Tstar shift, group bundle}%
\Index{texLie}
   {left structural constant of Lie algebra}%
   {left structural constant of Lie algebra}%
\Index{texVectorSpace}
   {left vector space}%
   {left vector space}%
\Index{texRepresentation}
   {left-side contravariant representation of group}%
   {left-side contravariant representation}%
\Index{texRepresentation}
   {left-side covariant representation of group}%
   {left-side covariant representation}%
\Index{texRepresentation}
   {left-side representation of group}%
   {left-side representation of group}%
\Index{texRepresentation}
   {left-side transformation}%
   {left-side transformation}%
\Index{texLie}
   {Lie algebra of Lie group}%
   {algebra Lie group Lie}%
\SubIndex{texLie}
   {left defined}%
   {left defined algebra Lie}%
\SubIndex{texLie}
   {right defined}%
   {right defined algebra Lie}%
\Index{texDiffProperty}
   {Lie derivative}%
   {Lie derivative}%
\SubIndex{texDiffProperty}
   {of connection}%
   {Lie derivative of connection}%
\SubIndex{texDiffProperty}
   {of metric}%
   {Lie derivative of metric}%
\Index{texLie}
   {Lie group basic operators}%
   {Lie group basic operators}%
\Index{texBundleRelation}
   {lift of correspondence}%
   {lift of correspondence}%
\Index{texAlgebraBundle}
   {lift of map}%
   {lift of map}%
\Index{texRepresentation}
   {linear representation of group}%
   {linear representation of group}%
\Index{texReferenceFrame}
   {linearly dependent vector fields}%
   {linearly dependent vector fields}%
\Index{texReferenceFrame}
   {local reference frame}%
   {local reference frame}%
\Index{texGenRelativity}
   {Lorentz transformation}%
   {Lorentz transformation}%
\SetIndexSpace%
\Index{texReferenceFrame}
   {map of type $G$ on manifold}%
   {map of type G on manifold}%
\Index{texVectorSpace}
   {matrix of \drc linear map}%
   {matrix of drc linear map}%
\Index{texAffine}
   {metric-affine manifold}%
   {metric-affine manifold}%
\Index{}
   {movement on basis manifold}%
   {movement transformation}%
\SetIndexSpace%
\Index{texPolymodule}
   {$(n)$\hyph vector space}%
   {(n)-vector space}%
\Index{texBundleRelation}
   {$n$\Hyph ary fibered relation}%
   {fibered relation}%
\Index{texAffine}
   {nonmetricity}%
   {nonmetricity}%
\Index{texVectorSpace}
   {nonsingular \drc system of linear equations}%
   {nonsingular system of linear equations}%
\Index{texRepresentation}
   {nonsingular transformation}%
   {nonsingular transformation}%
\Index{texCalculus}
   {norm on \drc vector space}%
   {norm on drc vector space}%
\Index{texCalculus}
   {normed \drc vector space}%
   {normed drc vector space}%
\SetIndexSpace%
\Index{texAlgebraBundle}
   {operation on bundle}%
   {operation on bundle}%
\Index{texBundleRelation}
   {opposite fibered preordering}%
   {opposite fibered preordering}%
\Index{texRepresentation}
   {orbit of representation of group}%
   {orbit of representation of group}%
\Index{texAlgebraBundle}
   {orbit of representation of group bundle}%
   {orbit of representation of group bundle}%
\Index{texDrcRepresentation}
   {orbit of \Ts representation of group}%
   {orbit of representation of group}%
\Index{texBasis}
   {orthonornal basis}%
   {Orthonornal Basis}%
\SetIndexSpace%
\Index{texReferenceFrame}
   {parallelogram}%
   {parallelogram}%
\Index{texCalculus}
   {partial derivative of mapping $f$ with respect to variable $v^{\gi a}$}%
   {partial derivative of mapping with respect to variable, skew field}%
\Index{texCalculus}
   {partial derivative of mapping $\Vector f$ with respect to variable $v^{\gi a}$}%
   {partial derivative of mapping with respect to variable, drc vector space}%
\Index{texBasis}
   {passive representation}%
   {passive representation}%
\Index{texBasis}
   {passive transformation on basis manifold}%
   {passive transformation}%
\Index{texDrcBasis}
   {passive \Ts representation}%
   {passive representation}%
\Index{texReferenceFrame}
   {pfaffian derivative}%
   {pfaffian derivative}%
\Index{texNewton}
   {potential energy}%
   {potential energy}%
\Index{texDrcBasis}
   {product of geometrical object and constant}%
   {product of geometrical object and constant}%
\Index{texBasis}
   {product of geometrical object and constant in vector space}%
   {product of geometrical object and constant, vector space}%
\SetIndexSpace%
\Index{texBasis}
   {quasi affine transformation on basis manifold}%
   {quasi affine transformation}%
\Index{texBasis}
   {quasi movement on basis manifold}%
   {quasi movement}%
\SetIndexSpace%
\Index{texBiring}
   {\sups row of matrix}%
   {r row}%
\Index{texVectorSpace}
   {\sups rows \drc vector space}%
   {sups rows drc vector space}%
\Index{texBiring}
   {$r$\hyph row of matrix}%
   {r-row}%
\Index{texBiring}
   {\RC inverse element of biring}%
   {rc-inverse element}%
\Index{texVectorSpace}
   {\RC major minor}%
   {RC-major minor}%
\Index{texVectorSpace}
   {\RC matrix group}%
   {rc-matrix group}%
\Index{texVectorSpace}
   {\RC nonsingular matrix}%
   {RC nonsingular matrix}%
\Index{texBiring}
   {\RC power}%
   {rc power}%
\Index{texBiring}
   {\RC product of matrices}%
   {rc-product of matrices}%
\Index{texBiring}
   {\RC quasideterminant}%
   {RC-quasideterminant}%
\Index{texVectorSpace}
   {\RC rank of matrix}%
   {rc-rank of matrix}%
\Index{texVectorSpace}
   {\RC singular matrix}%
   {RC singular matrix}%
\Index{texVectorSpace}
   {\rcd vector space}%
   {rcd vector space}%
\Index{texAlgebraBundle}
   {reduced Cartesian product of bundles}%
   {reduced Cartesian product of bundles}%
\Index{texAlgebraBundle}
   {reduced Cartesian product of total spaces}%
   {reduced Cartesian product of total spaces}%
\Index{texBundleRelation,texBundleRelation}
   {reduced fibered correspondence from $\mathcal{A}$ to $\mathcal{B}$}%
   {reduced fibered correspondence from A to B}%
\Index{texBundleRelation}
   {reduced fibered correspondence in $\mathcal{A}$}%
   {reduced fibered correspondence in A}%
\Index{texBiring}
   {reducible biring}%
   {reducible biring}%
\Index{texReferenceFrame}
   {reference frame}%
   {reference frame}%
\Index{texGenRelativity}
   {reference frame in event space}%
   {reference frame in event space}%
\Index{texReferenceFrame}
   {reference frame manifold}%
   {reference frame manifold}%
\Index{texBundleRelation}
   {reflexive $2$\Hyph ary fibered relation}%
   {reflexive 2 ary fibered relation}%
\Index{texDrcRepresentation}
   {representation of algebra $A$}%
   {}%
\SubIndex{texDrcRepresentation}
   {effective}%
   {effective representation of algebra}%
\SubIndex{texDrcRepresentation}
   {left-side}%
   {left-side representation of algebra}%
\SubIndex{texDrcRepresentation}
   {right-side}%
   {right-side representation of algebra}%
\SubIndex{texDrcRepresentation}
   {single transitive}%
   {single transitive representation of algebra}%
\SubIndex{texDrcRepresentation}
   {\sT}%
   {starT representation of algebra}%
\SubIndex{texDrcRepresentation}
   {transitive}%
   {transitive representation of algebra}%
\SubIndex{texDrcRepresentation}
   {\Ts}%
   {Tstar representation of algebra}%
\Index{texAlgebraBundle}
   {representation of algebra bundle}%
   {}%
\SubIndex{texAlgebraBundle}
   {effective}%
   {effective representation of algebra bundle}%
\SubIndex{texAlgebraBundle}
   {left-side}%
   {left-side representation of algebra bundle}%
\SubIndex{texAlgebraBundle}
   {right-side}%
   {right-side representation of algebra bundle}%
\SubIndex{texAlgebraBundle}
   {single transitive}%
   {single transitive representation of algebra bundle}%
\SubIndex{texAlgebraBundle}
   {\sT}%
   {starT representation of algebra bundle}%
\SubIndex{texAlgebraBundle}
   {transitive}%
   {transitive representation of algebra bundle}%
\SubIndex{texAlgebraBundle}
   {\Ts}%
   {Tstar representation of algebra bundle}%
\Index{texDrcRepresentation}
   {representation of group}%
   {}%
\SubIndex{texDrcRepresentation}
   {contravariant \Ts}%
   {Tstar contravariant representation of group}%
\SubIndex{texDrcRepresentation}
   {covariant \Ts}%
   {Tstar covariant representation of group}%
\SubIndex{texDrcBasis}
   {\drc linear \sT}%
   {linear representation of group}%
\SubIndex{texDrcRepresentation}
   {effective}%
   {effective representation of group}%
\SubIndex{texDrcBasis}
   {\rcd}%
   {rcd linear representation of group}%
\SubIndex{texDrcRepresentation}
   {\sT}%
   {starT representation of group}%
\SubIndex{texDrcRepresentation}
   {\Ts}%
   {Tstar representation of group}%
\Index{texRepresentation}
   {representation of group}%
   {representation of group}%
\Index{texAlgebraBundle}
   {representation of group bundle}%
   {}%
\SubIndex{texAlgebraBundle}
   {effective}%
   {effective representation of group bundle}%
\SubIndex{texAlgebraBundle}
   {\sT}%
   {starT representation of group bundle}%
\SubIndex{texAlgebraBundle}
   {\Ts}%
   {Tstar representation of group bundle}%
\SubIndex{texAlgebraBundle}
   {\Ts contravariant}%
   {Tstar contravariant representation of group bundle}%
\SubIndex{texAlgebraBundle}
   {\Ts covariant}%
   {Tstar covariant representation of group bundle}%
\Index{texDrcBasis}
   {representative of geometrical object in vector space}%
   {representative of geometrical object, vector space}%
\Index{texBasis}
   {representative of geometrical object in vector space}%
   {representative of geometrical object, vector space}%
\Index{texBundleRelation}
   {restriction of correspondence $\Phi$ to set $C$}%
   {restriction of correspondence}%
\Index{texLie}
   {right invariant vector field}%
   {right invariant vector}%
\Index{texVectorSpace}
   {right module}%
   {right module}%
\Index{texDrcRepresentation}
   {right shift}%
   {right shift}%
\Index{texRepresentation}
   {right shift on group}%
   {right shift, group}%
\Index{texLie}
   {right structural constant of Lie algebra}%
   {right structural constant of Lie algebra}%
\Index{texVectorSpace}
   {right vector space}%
   {right vector space}%
\Index{texRepresentation}
   {right-side contravariant representation of group}%
   {right-side contravariant representation}%
\Index{texRepresentation}
   {right-side covariant representation of group}%
   {right-side covariant representation}%
\Index{texRepresentation}
   {right-side representation of group}%
   {right-side representation of group}%
\Index{texRepresentation}
   {right-side transformation}%
   {right-side transformation}%
\Index{texAlgebraBundle}
   {ring bundle}%
   {ring bundle}%
\Index{texRepresentation}
   {row vector}%
   {row vector}%
\Index{texVectorSpace}
   {$R\star$\Hyph module}%
   {Rstar-module}%
\SetIndexSpace%
\Index{texNewton}
   {scalar potential}%
   {scalar potential}%
\Index{texNewton}
   {second Newton law}%
   {Second Newton law}%
\Index{texRepresentation}
   {single transitive representation of group}%
   {single transitive representation of group}%
\Index{texTidal}
   {speed of deviation}%
   {speed of deviation}%
\Index{texDrcBasis}
   {standard coordinates of basis}%
   {standard coordinates of basis}%
\Index{texBasis}
   {standard coordinates of basis}%
   {standard coordinates of basis}%
\Index{texBiring}
   {standard representation of matrix}%
   {Standard representation}%
\Index{texVectorSpace}
   {$\star D$\hyph vector space}%
   {starD-vector space}%
\Index{texLinearMap}
   {$\star D$\Hyph product of \drc linear map $A$ over scalar}%
   {starD product of drc linear map over scalar}%
\Index{texVectorSpace}
   {$\star R$\hyph module}%
   {starR-module}%
\Index{texDrcRepresentation}
   {\sT shift}%
   {starT shift}%
\Index{texAlgebraBundle}
   {\sT shift on group bundle}%
   {starT shift, group bundle}%
\Index{texAlgebraBundle}
   {subalgebra bundle}%
   {subalgebra bundle}%
\Index{texBundleRelation}
   {subbundle}%
   {subbundle}%
\Index{texLinearMap}
   {sum of \drc linear maps}%
   {sum of drc linear maps, drc vector spaces}%
\Index{texDrcBasis}
   {sum of geometrical objects}%
   {sum of geometrical objects}%
\Index{texBasis}
   {sum of geometrical objects in vector space}%
   {sum of geometrical objects, vector space}%
\Index{texBundleRelation}
   {symmetric $2$\Hyph ary fibered relation}%
   {symmetric 2 ary fibered relation}%
\Index{texBasis}
   {symmetry group}%
   {symmetry group}%
\Index{texDrcBasis}
   {symmetry group}%
   {SymmetryGroup}%
\Index{texGenRelativity}
   {synchronization of reference frame}%
   {synchronization of reference frame}%
\SetIndexSpace%
\Index{texLie}
   {tensor product of representations}%
   {tensor product of representations}%
\Index{texCalculus}
   {topological \drc vector space}%
   {topological drc vector space}%
\Index{texCalculus}
   {topological skew field}%
   {topological skew field}%
\Index{texAffine}
   {torsion}%
   {torsion}%
\Index{texDrcRepresentation}
   {transformation}%
   {}%
\SubIndex{texDrcRepresentation}
   {left-side}%
   {left-side transformation}%
\SubIndex{texDrcRepresentation}
   {nonsingular}%
   {nonsingular transformation}%
\SubIndex{texDrcRepresentation}
   {right-side}%
   {right-side transformation}%
\SubIndex{texDrcRepresentation}
   {\sT}%
   {starT transformation}%
\SubIndex{texDrcRepresentation}
   {\Ts}%
   {Tstar transformation}%
\Index{texDrcBasis}
   {transformation on basis manifold}%
   {}%
\SubIndex{texDrcBasis}
   {active}%
   {active transformation, vector space}%
\SubIndex{texTypeBasis}
   {affine}%
   {affine transformation}%
\SubIndex{texTypeBasis}
   {movement}%
   {movement transformation}%
\SubIndex{texDrcBasis}
   {passive}%
   {passive transformation, vector space}%
\SubIndex{texTypeBasis}
   {quasi affine}%
   {quasi affine transformation}%
\SubIndex{texTypeBasis}
   {quasi movement}%
   {quasi movement}%
\Index{texAlgebraBundle}
   {transformation on bundle}%
   {}%
\SubIndex{texAlgebraBundle}
   {left-side}%
   {left-side transformation of bundle}%
\SubIndex{texAlgebraBundle}
   {nonsingular}%
   {nonsingular transformation of bundle}%
\SubIndex{texAlgebraBundle}
   {\sT}%
   {starT transformation of bundle}%
\SubIndex{texAlgebraBundle}
   {\Ts}%
   {Tstar transformation of bundle}%
\Index{texBundleRelation}
   {transitive $2$\Hyph ary fibered relation}%
   {transitive 2 ary fibered relation}%
\Index{texRepresentation}
   {transitive representation of group}%
   {transitive representation of group}%
\Index{texVectorSpace}
   {\Ts linear composition of  vectors}%
   {linear composition of  vectors}%
\Index{texVectorSpace}
   {\Ts matrices vector space}%
   {matrices vector space}%
\Index{texDrcRepresentation}
   {\Ts shift}%
   {Tstar shift}%
\Index{texAlgebraBundle,texLinearMap,texDrcRepresentation}
   {twin representations}%
   {}%
\SubIndex{texDrcRepresentation}
   {of group}%
   {twin representations of group}%
\SubIndex{texAlgebraBundle}
   {of group bundle}%
   {twin representations of group bundle}%
\SubIndex{texLinearMap}
   {of skew field}%
   {twin representations of skew field}%
\Index{texReferenceFrame}
   {type $G$ reference frame}%
   {type G reference frame}%
\SetIndexSpace%
\Index{texVectorSpace}
   {unitarity law}%
   {}%
\SubIndex{texVectorSpace}
   {for \Ts vector space}%
   {unitarity law, Tstar vector space}%
\SetIndexSpace%
\Index{texPolymodule}
   {($D_1\RCstar$, ..., $D_n\RCstar$)\hyph vector space}%
   {(d1rc,dnrc)-vector space}%
\Index{texPolymodule}
   {($S\star$, $\star T$)\hyph vector space}%
   {(Sstar,starT)-vector space}%
\Index{texCalculus}
   {valued skew field}%
   {valued skew field}%
\Index{texAlgebraBundle}
   {vector bundle}%
   {vector bundle}%
\Index{texNewton}
   {vector potential}%
   {vector potential}%
\Index{texVectorSpace}
   {vector space type}%
   {vector space type}%

\CloseIndex

%% file: Symbol.English.tex
\def\indexname{Special Symbols and Notations}
\OpenIndex

\SetIndexSpace%
\Symb{texBiring}
   {$(^{\gi a}_{\gi b})$\hyph\CR quasideterminant}%
   {a b CR quasideterminant definition}%
\Symb{texBiring}
   {$(^{\gi a}_{\gi b})$\hyph \RC quasideterminant}%
   {a b RC-quasideterminant definition}%
\Symb{texBiring}
   {minor}%
   {A from b a}%
\Symb{texBiring}
   {minor}%
   {A from columns T}%
\Symb{texBiring}
   {minor}%
   {A from rows S}%
\Symb{texBiring}
   {minor}%
   {A without column a}%
\Symb{texBiring}
   {minor}%
   {A without columns T}%
\Symb{texBiring}
   {minor}%
   {A without row b}%
\Symb{texBiring}
   {minor}%
   {A without rows S}%
\Symb{texPolymodule}
   {active transformation}%
   {active transformation}%
\Symb{texTypeBasis}
   {affine space}%
   {affine space}%
\Symb{texBasis}
   {affine space}%
   {An}%
\Symb{texBiring}
   {\subs row ($c$\hyph row) of matrix}%
   {c row}%
\Symb{texBiring}
   {\CR power of element $A$ of biring}%
   {cr power}%
\Symb{texBiring}
   {\CR inverse element of biring}%
   {cr-inverse element}%
\Symb{texBiring}
   {\CR product of matrices}%
   {cr-product of matrices}%
\Symb{texVectorSpace}
   {\dcr vector}%
   {dcr vector}%
\Symb{texLie}
   {derivative of left shift}%
   {derivative of left shift}%
\Symb{texLie}
   {derivative of left shift}%
   {derivative of left shift, 1-Parameter Group}%
\Symb{texLie}
   {derivative of right shift}%
   {derivative of right shift}%
\Symb{texLie}
   {}%
   {derivative of right shift}%
\Symb{texLie}
   {derivative of right shift}%
   {derivative of right shift, 1-Parameter Group}%
\Symb{texLie}
   {derivative of left shift}%
   {derivative of Tstar shift}%
\Symb{texVectorSpace}
   {\drc vector}%
   {drc vector}%
\Symb{texAffine}
   {derivative}%
   {overline nabla_l, definition 2}%
\Symb{texPolymodule}
   {passive transformation}%
   {passive transformation}%
\Symb{texBiring}
   {\sups row ($r$\hyph row) of matrix}%
   {r row}%
\Symb{texBiring}
   {\RC power of element $A$ of biring}%
   {rc power}%
\Symb{texBiring}
   {\RC inverse element of biring}%
   {rc-inverse element}%
\Symb{texBiring}
   {\RC product of matrices}%
   {rc-product of matrices}%
\Symb{texBiring}
   {\RC quasideterminant}%
   {RC-quasideterminant definition}%
\Symb{texDrcRepresentation}
   {right shift}%
   {right shift}%
\Symb{texAlgebraBundle}
   {\sT shift}%
   {starT shift, group bundle}%
\Symb{texDrcRepresentation}
   {left shift}%
   {Tstar shift}%
\Symb{texAlgebraBundle}
   {\Ts shift}%
   {Tstar shift, group bundle}%
\Symb{texReferenceFrame}
   {anholonomic coordinates of vector}%
   {vector anholonomic coordinates}%
\Symb{texReferenceFrame}
   {holonomic coordinates of vector}%
   {vector holonomic coordinates}%

\SetIndexSpace%
\Symb{texBasis}
   {basis manifold of affine space}%
   {BAn}%
\Symb{texReferenceFrame}
   {basis manifold of manifold}%
   {basis manifold of manifold}%
\Symb{texPolymodule}
   {basis manifold of vector space}%
   {basis manifold of vector space}%
\Symb{texBasis}
   {basis manifold of vector space $\mathcal{V}$}%
   {basis manifold of vector space}%
\Symb{texBasis}
   {basis manifold of central affine space}%
   {BCAn}%
\Symb{texBasis}
   {basis manifold of Euclid space}%
   {BEn}%
\Symb{texTypeBasis}
   {basis manifold of affine space}%
   {FAn}%
\Symb{texTypeBasis}
   {basis manifold of central affine space}%
   {FCAn}%
\Symb{texTypeBasis}
   {basis manifold of Euclid space}%
   {FEn}%

\SetIndexSpace%
\Symb{texBasis}
   {central affine space}%
   {CAn}%
\Symb{texTypeBasis}
   {central affine space}%
   {central affine space}%
\Symb{texLie}
   {left structural constant of Lie algebra}%
   {left structural constant of Lie algebra}%
\Symb{texLie}
   {right structural constant of Lie algebra}%
   {right structural constant of Lie algebra}%
\Symb{texReferenceFrame}
   {set of functions defined on manifold}%
   {set of functions defined on manifold}%

\SetIndexSpace%
\Symb{texLieRepresentation}
   {basis vector of \sT representation}%
   {basis vector of starT representation}%
\Symb{texLieRepresentation}
   {basis vector of \sT representation}%
   {basis vector of starT representation, coordinates}%
\Symb{texLieRepresentation}
   {basis vector of \Ts representation}%
   {basis vector of Tstar representation}%
\Symb{texLieRepresentation}
   {basis vector of \Ts representation}%
   {basis vector of Tstar representation, coordinates}%
\Symb{texVectorSpace}
   {\subs rows \dcr vector space}%
   {c rows dcr vector space}%
\Symb{texCalculus}
   {differential of function}%
   {differential, drc vector space to drc vector space}%
\Symb{texCalculus}
   {differential of function}%
   {differential, drc vector space to skew field}%
\Symb{texVectorSpace}
   {\drc coordinate vector space}%
   {drc coordinate vector space}%
\Symb{texVectorSpace}
   {matrices vector space}%
   {matrices vector space}%
\Symb{texAffine}
   {Cartan derivative}%
   {overbrace D}%
\Symb{texAffine}
   {derivative}%
   {overline D}%
\Symb{texCalculus}
   {partial derivative of mapping $\Vector f$ with respect to variable $v^{\gi a}$}%
   {partial derivative of mapping, 1, drc vector space}%
\Symb{texCalculus}
   {partial derivative of mapping $f$ with respect to variable $v^{\gi a}$}%
   {partial derivative of mapping, 1, skew field}%
\Symb{texVectorSpace}
   {\sups rows \drc vector space}%
   {r rows drc vector space}%
\Symb{texTidal}
   {speed of deviation}%
   {speed of deviation}%
\Symb{texVectorSpace}
   {vector space type}%
   {vector space type}%

\SetIndexSpace%
\Symb{texTypeBasis}
   {affine basis}%
   {Affine Basis}%
\Symb{texBasis}
   {affine basis}%
   {Affine Basis}%
\Symb{texTypeBasis}
   {basis}%
   {basis}%
\Symb{texBasis}
   {basis of vector space}%
   {Basis e}%
\Symb{texBasis}
   {basis in vector space $\mathcal{V}$}%
   {basis in V}%
\Symb{texVectorSpace}
   {basis of vector space}%
   {basis, vector space}%
\Symb{texPolymodule}
   {basis of $(n)$\hyph vector space}%
   {basis,n vector space}%
\Symb{texAlgebraBundle}
   {Cartesian power of total spaces}%
   {Cartesian power of total spaces}%
\Symb{texAlgebraBundle}
   {Cartesian product of total spaces}%
   {Cartesian product of total spaces, definition 1}%
\Symb{texBasis}
   {central affine basis}%
   {Central Affine Basis}%
\Symb{texReferenceFrame}
   {form of reference frame}%
   {dual forms, reference frame}%
\Symb{texTypeBasis}
   {Euclid space}%
   {En}%
\Symb{texBasis}
   {Euclid space}%
   {En}%
\Symb{texBasis}
   {pseudo Euclid space}%
   {Enm}%
\Symb{texTypeBasis}
   {pseudo Euclid space}%
   {Enm}%
\Symb{texAlgebraBundle}
   {identical transformation of bundle}%
   {identical transformation of bundle}%
\Symb{texBasis}
   {orthonornal basis}%
   {Orthonornal Basis}%
\Symb{texAlgebraBundle}
   {reduced Cartesian product of total spaces}%
   {reduced Cartesian product of total spaces, definition 1}%
\Symb{texAlgebraBundle}
   {set of nonsingular \sT transformations of bundle $\mathcal{E}$}%
   {set of starT nonsingular transformations of bundle}%
\Symb{texAlgebraBundle}
   {set of nonsingular \Ts transformations of bundle $\mathcal{E}$}%
   {set of Tstar nonsingular transformations of bundle}%
\Symb{texBasis}
   {standard coordinates of basis}%
   {standard coordinates of basis}%
\Symb{texReferenceFrame}
   {standard coordinates of reference frame}%
   {standard coordinates of reference frame}%
\Symb{texReferenceFrame}
   {vector field of reference frame}%
   {vector field of reference frame}%
\Symb{texReferenceFrame}
   {vector field of reference frame}%
   {vector field, reference frame}%
\Symb{texBasis}
   {vector of basis}%
   {vector of basis}%

\SetIndexSpace%
\Symb{texVectorSpace}
   {coordinates of basis in \subs rows \dcr vector space}%
   {basis coordinates, c rows dcr vector space}%
\Symb{texVectorSpace}
   {coordinates of basis in \sups rows \drc vector space}%
   {basis coordinates, r rows drc vector space}%
\Symb{texVectorSpace}
   {basis for \subs rows \dcr vector space}%
   {basis, c rows dcr vector space}%
\Symb{texVectorSpace}
   {basis for \sups rows \drc vector space}%
   {basis, r rows drc vector space}%
\Symb{texDiffEq}
   {central affine basis}%
   {Central Affine Basis}%
\Symb{texReferenceFrame}
   {coordinate reference frame}%
   {coordinate reference frame}%
\Symb{texAlgebraBundle}
   {homomorphism of algebra bundles}%
   {homomorphism of algebra bundles}%
\Symb{texBundleRelation}
   {inverse fibered correspondence}%
   {inverse fibered correspondence, 1}%
\Symb{texBundleRelation}
   {inverse reduced fibered correspondence}%
   {inverse reduced fibered correspondence, 1}%
\Symb{texAlgebraBundle}
   {map to Cartesian product}%
   {map to Cartesian product}%
\Symb{texDrcRepresentation}
   {representation orbit of group $G$}%
   {orbit of representation of group}%
\Symb{texTypeBasis}
   {orthonornal basis}%
   {Orthonornal Basis}%
\Symb{texReferenceFrame}
   {reference frame}%
   {reference frame}%
\Symb{texReferenceFrame}
   {reference frame}%
   {reference frame, extensive definition}%
\Symb{texPolymodule}
   {standard coordinates of basis}%
   {standard coordinates of basis}%
\Symb{texPolymodule}
   {vector of basis}%
   {vector of basis}%

\SetIndexSpace%
\Symb{texVectorSpace}
   {\CR matrix group}%
   {cr-matrix group}%
\Symb{texLie}
   {algebra Lie of group Lie}%
   {g}%
\Symb{texLie}
   {left defined algebra Lie of group Lie}%
   {gl}%
\Symb{texTypeBasis}
   {affine transformation group}%
   {GLAn}%
\Symb{texBasis}
   {affine transformation group}%
   {GLAn}%
\Symb{texLie}
   {right defined algebra Lie of group Lie}%
   {gr}%
\Symb{texBasis}
   {group of homomorphisms of vector space $\mathcal{V}$}%
   {GV}%
\Symb{texDrcRepresentation}
   {orbit of effective covariant \sT representation of group}%
   {orbit of effective starT covariant representation of group}%
\Symb{texAlgebraBundle}
   {orbit of effective covariant \sT representation of group bundle}%
   {orbit of effective starT covariant representation of group bundle}%
\Symb{texDrcRepresentation}
   {orbit of effective covariant \Ts representation of group}%
   {orbit of effective Tstar covariant representation of group}%
\Symb{texAlgebraBundle}
   {orbit of effective covariant \Ts representation of group bundle}%
   {orbit of effective Tstar covariant representation of group bundle}%
\Symb{texVectorSpace}
   {\RC matrix group}%
   {rc-matrix group}%

\SetIndexSpace%
\Symb{texBiring}
   {Hadamard inverse of matrix}%
   {Hadamard inverse of matrix}%
\Symb{texLinearMap}
   {\rcd vector space of \drc linear maps}%
   {rcd vector space of drc linear maps}%

\SetIndexSpace%
\Symb{texLieRepresentation}
   {infinitesimal generator of representation}%
   {infinitesimal generator of representation}%
\Symb{texLinearLie}
   {Lie group infinitesimal generators}%
   {Lie group infinitesimal generators}%

\SetIndexSpace%
\Symb{texRepresentation}
   {left shift}%
   {left shift}%
\Symb{texDiffProperty}
   {Lie derivative of connection}%
   {Lie derivative of connection}%
\Symb{texDiffProperty}
   {Lie derivative of metric}%
   {Lie derivative of metric}%
\Symb{texBasis}
   {passive transformation}%
   {passive transformation}%
\Symb{texRepresentation}
   {set of left-side nonsingular transformations of set $M$}%
   {set of left-side nonsingular transformations}%

\SetIndexSpace%
\Symb{texDrcRepresentation}
   {set of \sT nonsingular \sT transformations of set $M$}%
   {set of starT nonsingular transformations}%
\Symb{texDrcRepresentation}
   {set of nonsingular \Ts transformations of set $M$}%
   {set of Tstar nonsingular transformations}%

\SetIndexSpace%
\Symb{texBasis}
   {geometrical object in coordinate representation}%
   {geometrical object, coordinate vector space}%
\Symb{texBasis}
   {geometrical object}%
   {geometrical object, vector space}%
\Symb{texRepresentation}
   {orbit of representation of the group $G$}%
   {orbit of representation of group}%
\Symb{texAlgebraBundle}
   {orbit of representation of group bundle $\mathcal{G}$}%
   {orbit of representation of group bundle}%

\SetIndexSpace%
\Symb{texAlgebraBundle}
   {bundle}%
   {bundle}%
\Symb{texAlgebraBundle}
   {Cartesian power of bundle}%
   {Cartesian power of bundle}%
\Symb{texAlgebraBundle}
   {Cartesian product of bundles}%
   {Cartesian product of bundles, definition 1}%
\Symb{texAlgebraBundle}
   {reduced Cartesian product of bundles}%
   {reduced Cartesian product of bundles, definition 1}%
\Symb{texAlgebraBundle}
   {set of nonsingular \sT transformations of bundle $\bundle{}pE{}$}%
   {set of starT nonsingular transformations of bundle, projection}%
\Symb{texAlgebraBundle}
   {set of nonsingular \Ts transformations of bundle $\bundle{}pE{}$}%
   {set of Tstar nonsingular transformations of bundle, projection}%

\SetIndexSpace%
\Symb{texBasis}
   {active transformation}%
   {active transformation}%
\Symb{texAffine}
   {Cartan curvature}%
   {Cartan curvature}%
\Symb{texVectorSpace}
   {\CR rank of matrix}%
   {cr-rank of matrix}%
\Symb{texBundleRelation}
   {diagonal in bundle  $\bundle{}pA{}$}%
   {diagonal in bundle, 2}%
\Symb{texBundleRelation}
   {diagonal in bundle $\mathcal{A}$}%
   {diagonal in reduced bundle, 2}%
\Symb{texAffine}
   {curvature}%
   {GLn curvature_overline}%
\Symb{texVectorSpace}
   {\RC rank of matrix}%
   {rc-rank of matrix}%
\Symb{texRepresentation}
   {right shift}%
   {right shift}%
\Symb{texRepresentation}
   {set of right-side nonsingular transformations of set $M$}%
   {set of right-side nonsingular transformations}%

\SetIndexSpace%
\Symb{texBundleRelation}
   {composition of fibered correspondences}%
   {composition of fibered correspondences}%
\Symb{texBundleRelation}
   {inverse fibered correspondence}%
   {inverse fibered correspondence, 2}%
\Symb{texBundleRelation}
   {inverse reduced fibered correspondence}%
   {inverse reduced fibered correspondence, 2}%
\Symb{texVectorSpace}
   {linear span in vector space}%
   {linear span, vector space}%

\SetIndexSpace%
\Symb{texLie}
   {tangent plane to group $G$}%
   {TaG}%

\SetIndexSpace%
\Symb{texBasis}
   {coordinate vector space}%
   {coordinate vector space}%
\Symb{texBasis}
   {coordinates in vector space}%
   {coordinates in vector space}%
\Symb{texVectorSpace}
   {\dcr vector space}%
   {left CR vector space}%
\Symb{texVectorSpace}
   {\drc vector space}%
   {left RC vector space}%
\Symb{texLinearMap}
   {($S$, $T$)\hyph bimodule}%
   {R S bimodule}%
\Symb{texVectorSpace}
   {\crd vector space}%
   {right CR vector space}%
\Symb{texVectorSpace}
   {\rcd vector space}%
   {right RC vector space}%
\Symb{texBasis}
   {vector space}%
   {V}%
\Symb{texReferenceFrame}
   {vector space of vector fields}%
   {vector space of vector fields}%

\SetIndexSpace%
\Symb{texPolymodule}
   {geometrical object in coordinate representation		defined in vector space}%
   {geometrical object, coordinate vector space}%
\Symb{texPolymodule}
   {geometrical object in vector space}%
   {geometrical object, vector space}%

\SetIndexSpace%
\Symb{texReferenceFrame}
   {anholonomic coordinate}%
   {x(k)}%

\SetIndexSpace%
\Symb{texBundleRelation}
   {diagonal in bundle $\mathcal{A}$}%
   {diagonal in bundle, 1}%

\SetIndexSpace%
\Symb{texTidal}
   {deviation of trajectories}%
   {deviation of trajectories}%
\Symb{texDrcRepresentation}
   {identical transformation}%
   {identical transformation}%
\Symb{texRepresentation}
   {identical transformation}%
   {identical transformation}%
\Symb{texBasis}
   {image of vector $\Vector e_k\in\Basis e$ under isomorphism to coordinate vector space}%
   {image of vector e_k, coordinate vector space}%
\Symb{texBiring}
   {Kronecker symbol}%
   {Kronecker symbol}%

\SetIndexSpace%
\Symb{texReferenceFrame}
   {anholonomic coordinates of connection}%
   {anholonomic coordinates of connection}%
\Symb{texAffine}
   {Cartan symbol}%
   {Cartan symbol}%
\Symb{texAffine}
   {connection}%
   {conection overline}%
\Symb{texAffine}
   {Cartan connection}%
   {overbrace Gamma i kl}%
\Symb{texAlgebraBundle}
   {set of sections of bundle}%
   {set of sections of bundle}%

\SetIndexSpace%
\Symb{texLie}
   {inverse operator to operator $\psi_l$}%
   {inverse operator to operator psi l}%
\Symb{texLie}
   {inverse operator to operator $\psi_r$}%
   {inverse operator to operator psi r}%

\SetIndexSpace%
\Symb{texReferenceFrame}
   {anholonomity object}%
   {anholonomity object}%

\SetIndexSpace%
\Symb{texLie}
   {basic operator of group Lie}%
   {Lie Basic Operator L}%
\Symb{texLie}
   {}%
   {Lie Basic Operator L}%
\Symb{texLie}
   {basic operator of group Lie}%
   {Lie Basic Operator L, 1-Parameter Group}%
\Symb{texLie}
   {basic operator of group Lie}%
   {Lie Basic Operator R}%
\Symb{texLie}
   {}%
   {Lie Basic Operator R}%
\Symb{texLie}
   {basic operator of group Lie}%
   {Lie Basic Operator R, 1-Parameter Group}%

\SetIndexSpace%
\Symb{texReferenceFrame}
   {coordinate reference frame}%
   {coordinate reference frame, extensive definition}%
\Symb{texCalculus}
   {partial derivative of mapping $\Vector f$ with respect to variable $v^{\gi a}$}%
   {partial derivative of mapping, 2, drc vector space}%
\Symb{texCalculus}
   {partial derivative of mapping $f$ with respect to variable $v^{\gi a}$}%
   {partial derivative of mapping, 2, skew field}%
\Symb{texReferenceFrame}
   {derivative $e_{(k)}$}%
   {partial(k)}%

\SetIndexSpace%
\Symb{texLie}
   {Lie group composition law}%
   {Lie group composition law}%

\SetIndexSpace%
\Symb{texAffine}
   {Cartan derivative}%
   {overbrace nabla_l}%
\Symb{texAffine}
   {}%
   {overbrace nabla_l}%
\Symb{texAffine}
   {derivative}%
   {overline nabla_l, definition 1}%

\SetIndexSpace%
\Symb{texBundleRelation}
   {restriction of correspondence $\Phi$ to set $C$}%
   {restriction of correspondence}%

\SetIndexSpace%
\Symb{texAlgebraBundle}
   {Cartesian product of bundles}%
   {Cartesian product of bundles, definition 2}%
\Symb{texAlgebraBundle}
   {Cartesian product of total spaces}%
   {Cartesian product of total spaces, definition 2}%
\Symb{texAlgebraBundle}
   {reduced Cartesian product of bundles}%
   {reduced Cartesian product of bundles, definition 2}%
\Symb{texAlgebraBundle}
   {reduced Cartesian product of total spaces}%
   {reduced Cartesian product of total spaces, definition 2}%

\SetIndexSpace%
\Symb{texBundleRelation}
   {fibered subset}%
   {fibered subset}%
\Symb{texBundleRelation}
   {subbundle}%
   {subbundle}%

\CloseIndex